\documentclass[a4paper]{article}
\usepackage{geometry}
 \usepackage{ulem}

\usepackage{times,cite}
\usepackage[T1]{fontenc}
\usepackage[utf8]{inputenc}
\usepackage{graphicx}
\usepackage[dvipsnames]{xcolor}
\usepackage{amsmath,amsfonts,mathtools, amsthm, amssymb}
\usepackage{dsfont,xfrac,sidecap,caption}
\usepackage{nicefrac}

\usepackage[ruled, vlined]{algorithm2e}

\usepackage{cleveref}
\usepackage{todonotes}
\newtheorem{lemma}{Lemma}
\newtheorem{theorem}{Theorem}

\usepackage{multirow}

\usepackage{pgfplots}

\usepackage{placeins}

\newcommand{\nt}{\vec{n}}

\newcommand{\vt}{\mathbf{v}}

\newcommand{\ut}{\mathbf{u}}
\newcommand{\wt}{\mathbf{w}}
\newcommand{\sigmat}{\boldsymbol{\sigma}}
\newcommand{\Sigmat}{\boldsymbol{\Sigma}}
\newcommand{\Et}{\mathbf{E}}

\newcommand{\Ft}{\mathbf{F}}

\newcommand{\bs}[1]{\boldsymbol{#1}}

\renewcommand{\div}{{\operatorname{div}}}

\newcommand{\FL}{{\cal F}}
\newcommand{\SO}{{\cal S}}
\newcommand{\IN}{\Gamma}

\newtheorem{varform}{Variational formulation}
\newtheorem{remark}{Remark}

\newcommand{\sfrei}{}
\newcommand{\refone}{}
\newcommand{\reftwo}{}

\definecolor{dred}{rgb}{0.8,0,0}

\begin{document}

\title{Towards parallel time-stepping for the numerical simulation of atherosclerotic plaque growth}

\date{}

\author{Stefan Frei\thanks{Department of Mathematics and Statistics, University of Konstanz, 78457 Konstanz, Germany, stefan.frei@uni-konstanz.de (Corresponding author)},\, Alexander Heinlein\thanks{Delft University of Technology, Delft Institute of Applied Mathematics, Mekelweg 4, 2628CD Delft, Netherlands, a.heinlein@tudelft.nl}}

\maketitle                   

\begin{abstract}
The numerical simulation of atherosclerotic plaque growth is computationally prohibitive, since it involves a complex cardiovascular fluid-structure interaction (FSI) problem with a characteristic time scale of milliseconds to seconds, as well as a plaque growth process governed by reaction-diffusion equations, which takes place over several months.
In this work we combine a temporal homogenization approach, which separates the problem in computationally expensive FSI problems on a micro scale and a reaction-diffusion problem on the macro scale, with parallel time-stepping algorithms.
It has been found in the literature that parallel time-stepping algorithms do not perform well when applied directly to the FSI problem. To circumvent this problem, a parareal algorithm is applied on the macro-scale reaction-diffusion problem instead of the micro-scale FSI problem.
We investigate modifications in the coarse propagator of the parareal algorithm, in order to further reduce the number of costly micro problems to be solved. The approaches are tested in detailed numerical investigations based on serial simulations.
\end{abstract}

\section{Introduction}

Cardiovascular diseases are by far the most common cause of death in industrialized nations today. 
One of the most common cardiovascular diseases is the growth of plaque (atherogenesis) in coronary arteries or the pathological accumulation of plaque (atherosclerosis), which can result in heart attacks or strokes that are often fatal. Since the formation of plaque ranges over a long period of time, from months to several years, early diagnosis and treatment to prevent plaque growth can have a good chance of success. 

An important driving force
for plaque growth is the wall shear stress distribution, which varies significantly within each heartbeat, i.e., every second; see, for instance, \cite{Wahle:2006:PDV,Dhawan:2010:SSP} for a discussion on the dependence of plaque growth on the wall shear stress. Using three-dimensional fluid-structure interaction (FSI) simulations with realistic material models, the heartbeat scale has to be resolved by time steps in the order of milliseconds. We refer to~\cite{balzani2014aspects,Balzani:2016:NMF} for the first studies on using complex wall models that take into account the effect of the reinforcing fibers within three-dimensional FSI simulations. For further references on FSI for cardiovascular applications, see, e.g.,~\cite{formaggia2010cardiovascular,barker2010scalable,bazilevs2013computational,Hron:2006:MFE,kuttler2010coupling}; for FSI simulations in the context of plaque growth, see also~\cite{Figueroaetal2009, YangRichterJaegerNeussRadu2015, FreiRichterWick2016, ThonHemmleretal2018, PozziRedaelliVergara2021, Frei:2021:THN}.

Hence, realistic numerical simulations of plaque growth over several months that resolve each individual heartbeat would easily require $\mathcal{O}(10^9)$ sequential time steps. Consequently, even using today's fastest supercomputers, such a simulation is clearly infeasible. 

Therefore, in~\cite{FreiRichter2020} a temporal homogenization approach for separating the plaque growth time scale (macro scale) in the order of days and the FSI time scale (micro scale) in the order of milliseconds has been introduced. The approach is based on the assumption that the FSI is approximately periodic in time 
on the micro scale, which makes it possible to upscale the fluid dynamics to the macro scale. 
Using this approach it is possible to simulate only a moderate number of FSI time steps for each time step of the plaque growth problem. This means that we have to simulate only a few heartbeats (seconds) of the full FSI problem instead of the whole 24 hours of each day; the total number of time steps is reduced by a factor of $\mathcal{O}(10^4)$. Note that, completely neglecting the fine scale may introduce a significant error; see~\cite{FreiRichterWick2016,Frei:2021:THN}. See also the recent PhD thesis of Florian Sonner for further details on temporal homogenization for plaque growth~\cite{SonnerPhD}.

However, considering realistic three-dimensional simulations the number of time steps corresponding to fine scale FSI problems remains still infeasibly high, even after temporal homogenization. 
In order to further reduce the computational times, we introduce a new approach based on parallel time-stepping. We focus on a classical parallel time-stepping method, the parareal algorithm, which has been introduced by Lions et al.\ in~\cite{Lions:2001:PIT}{\reftwo; for an overview over the vast literature on parallel-in-time integration methods, we refer to the review papers~\cite{Gander:2015:50Y,ong_applications_2020} and the references therein.}

Due to a phase-shift in a coarse solution for hyperbolic partial differential equations (PDEs), such as the structural problem in FSI, it is generally challenging to apply parallel time-stepping methods to FSI problems; see, e.g.,~\cite{Margenberg:2021:PTS,ruprecht2018wave}.
Thus, instead of applying the parareal algorithm to the micro-scale FSI problem, we apply it to the homogenized plaque growth problem on the macro scale. Plaque growth is typically modeled by a system of reaction-diffusion equations; see, e.g.,~\cite{YangJaegerNeussRaduRichter2015,SilvaJaegerNeussRaduSequeira2020}. Instead of being hyperbolic, they have a parabolic character and hence can be solved more efficiently by parallel time-stepping methods. 

The focus of this paper is not the computation of accurate plaque-growth predictions in patient-specific geometries. Instead, our objective is to introduce a numerical framework for making the simulation times feasible and to investigate the methods numerically. Therefore, we make several simplifications: We focus on two-dimensional FSI simulations on a simple geometry. Furthermore, we consider two simplified models for the plaque growth, that is, the ordinary differential equation (ODE) model already considered in~\cite{FreiRichter2020,Frei:2021:THN} as well as a more complex partial differential equation (PDE) of reaction-diffusion type.
We formulate 
the parareal algorithm for the two-scale problem of plaque growth and also propose variants
to increase the efficiency of the coarse-scale propagators. 
Finally, we investigate the parallel time-stepping methods numerically based on serial simulations. As, even in two dimensions, the micro-scale FSI problems are generally much more expensive than the coarse plaque growth problem, we are able to give some good estimates for the speedup that can be expected in fully parallel simulations. 

The paper is structured as follows: In~\cref{sec:equations}, we introduce our fluid-structure interaction problem as well as the two solid growth models considered for modeling the plaque growth: an ODE-based model in~\cref{sec:ODE} and a reaction-diffusion equation-based model in~\cref{sec:PDE}. Next, in~\cref{sec:num}, we briefly introduce the variational formulation of the FSI problem and then describe the temporal homogenization approach including some first numerical results. In~\cref{sec:parallel}, we recap the parareal algorithm and discuss how we estimate the computational costs and the possible speed up in parallel simulations. We also discuss numerical results for the parallel time-stepping approach using the simple ODE plaque growth model. Furthermore, we introduce some ideas for reducing the costs of the coarse-scale propagators. In~\cref{sec:reaction-diffusion}, we show results for the reaction-diffusion plaque growth model, including a theoretical discussion on the computational costs of the algorithms. We conclude in~\cref{sec:conclusion}, where we also give a brief outlook on topics for future work.

\section{Model equations} \label{sec:equations}

We consider a time-dependent fluid-structure interaction system, where the fluid is modeled by the Navier--Stokes equations and the solid by the Saint Venant--Kirchhoff model. In order to account for the solid growth, we add a multiplicative growth term to the deformation gradient, which is motivated by typical plaque growth models~\cite{YangJaegerNeussRaduRichter2015,FreiRichterWick2016,MizerskiRichter2019, SilvaJaegerNeussRaduSequeira2020}.

\subsection{Fluid-structure interaction}

Consider a partition of an overall domain $\Omega(t) = \FL(t) \cup \IN(t)\cup \SO(t)$ into a fluid part $\FL(t)$, an interface $\IN(t)$ and a solid part $\SO(t)$. The blood flow and its interaction with the surrounding vessel wall is modelled by the following FSI system:
{\begin{equation}\label{fullplaquemodel}
  \begin{aligned}
      \rho_f(\partial_t \vt_f+\vt_f\cdot\nabla\vt_f) -
      \operatorname{div}\,\sigmat_f&= 0, &\qquad
      \operatorname{div}\,\vt_f&=0 \qquad\text{ in }\FL(t),\\
      \rho_s \partial_t\hat{\vt}_s -
      \operatorname{div}\,(\hat{\Ft}_e\hat{\Sigmat}_e) &=0, &\qquad
      \partial_t \hat\ut_s - \hat\vt_s &=0 \qquad \text{ in }\hat{S},\\
      \sigmat_f\nt_f+\sigmat_s\nt_s &=0, &\qquad
      \vt_f&=\vt_s \quad\;\; \text{ on }\IN(t).
  \end{aligned}
\end{equation}}
Here, $\vt_f$ and $\hat{\vt}_s$ denote the fluid and solid velocity, respectively, and 
$\hat{\ut}_s$ the solid displacement.
Quantities with a ``hat'' are defined in Lagrangian coordinates, while quantities without a ``hat'' are defined in the current Eulerian coordinate framework. Two quantities $\hat{f}(\hat{x})$ and $f(x)$ correspond to each other by a $C^{1,1}$-diffeomorphism $\hat{\xi}: \hat{\Omega}\to \Omega(t)$ and the relation $\hat{f} = f\circ \hat{\xi}$. Later, we will also need the solid deformation gradient $\hat{F}_s = I + \hat{\nabla } \hat{u}_s$, which is the derivative of $\hat{\xi}$ in the solid part.
The constants $\rho_f$ and $\rho_s$ are the densities of blood and vessel
wall, and $\nt_f$ and $\nt_s$ are outward pointing normal vectors of the fluid and solid 
domain, respectively.

By $\sigmat_f$ and $\sigmat_s$ we denote the Cauchy stress tensors of fluid and solid. Using the well-known Piola transformation between the Eulerian and the Lagrangian coordinate systems, we can relate $\sigmat_s$ and the second Piola--Kirchhoff stress $\hat{\Sigma}_e$ as follows: 
\begin{align*}
 \sigmat_s(x) = \hat{\sigmat}_s(\hat{x}) = \hat{J}_e^{-1} \hat{\Ft}_e\hat{\Sigmat}_e(\hat{x}) \hat{\Ft}_e^T, 
\end{align*}
where $\hat{F}_e$ is the elastic part of the deformation gradient $\hat{F}_s$,
For modeling the material behavior of the vessel wall, different approaches have been proposed and investigated in literature; see, e.g.,~\cite{Balzani:2016:NMF,Balzani:2006:PFS,brands2008modelling,Balzani:CAW:2022,holzapfel2000new} for more sophisticated material models, for instance, incorporating an anisotropic behavior due to the reinforcing fibers.
For the sake of simplicity, we use in this work the relatively simple Saint Venant-Kirchhoff model with the Lam\'e
material parameters $\mu_s$ and $\lambda_s$
\begin{equation}\label{structure-stress-lagrange}
  \hat \Sigmat_e= 2\mu_s \hat \Et_e
  +\lambda_s \operatorname{tr}(\hat\Et_e)I,\quad 
  \hat\Et_e:=\frac{1}{2}(\hat\Ft_e^T\hat\Ft_e-I).
\end{equation}
The Saint Venant-Kirchhoff model is based on Hooke's linear material law in a large strain formulation, resulting in a weakly nonlinear material model. 

\begin{figure}
\centering
\begin{tikzpicture}
\draw[fill=gray!30,line width=1pt] (-5,-2) -- (5,-2) -- (5,-1) -- (-5,-1) -- cycle;
\node () at (2.5,-1.5) {$\mathcal{S}$};
\draw[line width=1pt] (-5,-1) -- (5,-1) -- (5,1) -- (-5,1) -- cycle;
\draw[fill=gray!30,line width=1pt] (-5,1) -- (5,1) -- (5,2) -- (-5,2) -- cycle;
\node () at (2.5,1.5) {$\mathcal{S}$};

\node () at (3.0,0.5) {$\mathcal{F}$};

\draw[|<->|,line width=0.5pt] (-3.5,-2) -- (-3.5,-1); \node () at (-3.75,-1.5) {1};
\draw[|<->|,line width=0.5pt] (-3.5,-1) -- (-3.5,0); \node () at (-3.75,-0.5) {1};
\draw[|<->|,line width=0.5pt] (-3.5,0) -- (-3.5,1); \node () at (-3.75,0.5) {1};
\draw[|<->|,line width=0.5pt] (-3.5,1) -- (-3.5,2); \node () at (-3.75,1.5) {1};

\draw[|<->|,line width=0.5pt] (-5,-2.25) -- (0,-2.25); \node () at (-2.5,-2.5) {5};
\draw[|<->|,line width=0.5pt] (0,-2.25) -- (5,-2.25); \node () at (2.5,-2.5) {5};

\node () at (-5.5,1.5) {$\Gamma_s$};
\node () at (5.5,1.5) {$\Gamma_s$};
\node () at (-5.5,-1.5) {$\Gamma_s$};
\node () at (5.5,-1.5) {$\Gamma_s$};

\node () at (-5.5,0.25) {$\Gamma_f^{\text{in}}$};
\node () at (5.5,0.25) {$\Gamma_f^{\text{out}}$};

\node () at (-2.0,0.75) {$\Gamma$};
\node () at (-2.0,-0.75) {$\Gamma$};

\draw[line width=1.0pt,dotted] (-5.5,0) -- (5.5,0);
\draw[line width=1.0pt,dotted] (0,-2.5) -- (0,2.5);

\draw[fill=black] (0,-1) circle (0.1); \node () at (-0.5,0.75) {(0,1)};
\draw[fill=black] (0,1) circle (0.1); \node () at (-0.5,-0.75) {(0,-1)};

\draw[fill=black] (0,0) circle (0.05); \node () at (0.5,0.25) {(0,0)};
\end{tikzpicture} 
\caption{Sketch of the computational domain centered at the origin $(0,0)$; $\FL$ and $\SO$ are the fluid and solid parts, respectively, and the solid lines correspond to the fluid-solid interface $\Gamma$. The plaque growth is initiated at $(0,\pm1)$.
\label{fig:conf}}
\end{figure}
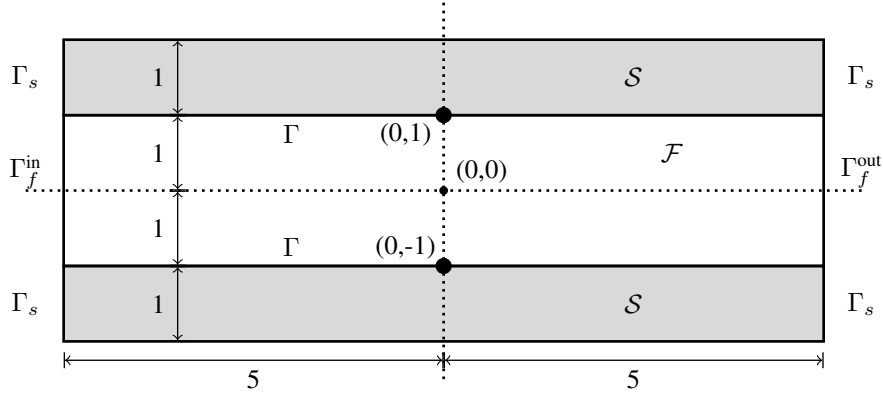

The blood flow is modeled as an incompressible Newtonian fluid, such that the
Cauchy stresses are given by
\begin{equation}\label{fluidstresses-euler}
  \sigmat_f=\rho_f\nu_f(\nabla\vt_f+\nabla\vt_f^T)-p_f I,
\end{equation}
where $\nu_f$ is the kinematic viscosity of blood. A sketch of the computational domain is given in~\cref{fig:conf}. We split the outer boundary of $\Omega$ into a solid part $\Gamma_s$ with homogeneous Dirichlet conditions, a fluid part $\Gamma^\text{in}_f$ with an inflow Dirichlet condition and an outflow part $\Gamma^\text{out}_f$, where a \textit{do-nothing condition} is imposed.

The boundary data is given by 
  \begin{align}\label{problem:long-scale-boundary}
    \vt_f = \vt^\text{in} \,\text{ on }\Gamma^\text{in}_f,
    \quad
    \rho_f\nu_f(\nt_f\cdot\nabla)\vt_f-p_f\nt_f =0
    & \,\text{ on }\Gamma^\text{out}_f,\quad
    \hat{\ut}_s=0\,\text{ on }\hat{\Gamma}_s,
  \end{align}
  where $\vt^\text{in}$ is the inflow velocity on $\Gamma^\text{in}_f$.
  
\subsection{Modelling of solid growth} \label{sec:solid_growth}

Developing a realistic model of plaque growth at the vessel wall is a complex task that involves the interaction of many different molecules 
and species, see for example~\cite{SilvaJaegerNeussRaduSequeira2020}. Furthermore, the plaque growth will also strongly depend on the geometry and material model of the arterial wall.
In this contribution, our focus does not lie on a realistic modeling of plaque growth, and thus, our model is greatly simplified. However, the two models considered here are chosen to behave similarly (from a numerical viewpoint) to more sophisticated models. We consider a simple ODE-based model in~\cref{sec:ODE} and a more complex but still relatively simple PDE model of reaction-diffusion type in~\cref{sec:PDE}. Both models focus on the influence of the concentration of foam cells $c_s$ on the growth. Numerical results can be found in~\cref{sec:numex1,sec:pararealres,sec:pararealmacro} as well as in \cref{sec:reaction-diffusion}, respectively.

\subsubsection{A simple ODE model for solid growth}
\label{sec:ODE}

In the first model, which is taken from~\cite{FreiRichterWick2016, Frei:2021:THN}, the evolution of the foam cell concentration $c_s$ depends only on its current value but not on its spatial distribution. Hence, it can be described by a simple ODE. The rate of formation of these cells depends on the distribution of the wall shear stress $\sigmat_f^{WS}$ at the vessel wall. As a result, we obtain the following simplified ODE model:
\begin{equation}\label{reaction}
\begin{aligned}
    \partial_t c_s & = \gamma(\sigmat_f^{WS},c_s)\coloneqq \alpha
    \big(1+c_s\big)^{-1}\left(1+\frac{\left\|\sigmat_f^{WS}\right\|_{L^2(\Gamma)}^2}{\sigma_0^{2}}\right)^{-1}, \\ 
    \sigmat_f^{WS} & \coloneqq 
    \rho\nu
    \big(I_d-\nt_f\nt_f^T\big)(\nabla\vt+\nabla\vt^T)\nt_f.
\end{aligned}
\end{equation}
The reference wall shear stress $\sigma_0$ and the scale separation parameter $\alpha$ are parameters of the growth model. For cardiovascular plaque growth, we have typically $\alpha = \mathcal{O}(10^{-7})\,\text{s}^{-1}$; see also~\cite{FreiRichter2020}. 

We model the solid growth by a multiplicative splitting of the 
 deformation gradient $\hat \Ft_s$
into an elastic part $\hat\Ft_e$ and a growth function $\hat\Ft_g$
 \begin{equation}\label{elasticgrowth-reference}
  \hat \Ft_s=\hat\Ft_e\hat\Ft_g\quad\Leftrightarrow\quad
  \hat\Ft_e = \hat\Ft_s\hat\Ft_g^{-1} =
  [I+\hat\nabla\hat\ut_s]\hat\Ft_g^{-1};
\end{equation}
cf.~\cite{RodriguezHogerMcCulloch1994, YangJaegerNeussRaduRichter2015, FreiRichterWick2016}.
In the ODE model, we use the following growth function depending on $c_s$
\begin{equation}\label{num:1:growth}
  \hat g(\hat x,\hat y,t) =  1+ c_s \exp\left(-\hat
  x^2\right)(2-|\hat y|),\quad 
  \hat\Ft_g(\hat x,\hat y,t):=\hat g(\hat x,\hat y,t)\,I.
\end{equation}
This means that the shape and position of the plaque growth is prescribed, but the growth rate depends on the variable $c_s$. As the simulation domain is centered around the origin, see~\cref{fig:conf}, growth is concentrated at $(0,\pm1)$, 
in the center of the domain.
It follows that
\begin{equation}
  \hat \Ft_g = \hat gI\quad\Rightarrow\quad
  \hat\Ft_e:=\hat g^{-1} \hat\Ft_s,\label{FgFe}
\end{equation}
and the elastic Green--Lagrange strain is given by
\begin{equation}\label{Es}
  \hat{\bs{E}}_e = \frac{1}{2}(\hat
      {\bs{F}}_e^T\hat{\bs{F}}_e-I) = \frac{1}{2}(\hat g^{-2}\hat
      {\bs{F}}_s^T\hat{\bs{F}}_s-I)
\end{equation}
resulting in the Piola--Kirchhoff stresses
\begin{equation}\label{Piola}
  \hat\Ft_e\hat\Sigmat_e = 2\mu_s \hat{\bs{F}}_e\hat{\bs{E}}_e +
  \lambda_s  \operatorname{tr}(\hat{\bs{E}}_e)\hat{\bs{F}}_e
  = 2\mu_s \hat g^{-1}\hat{\bs{F}}_s\hat{\bs{E}}_e +
  \lambda_s  \hat g^{-1}\operatorname{tr}(\hat{\bs{E}}_e)\hat{\bs{F}}_s.
\end{equation}

\subsubsection{A PDE reaction-diffusion model}
\label{sec:PDE}

Secondly, we consider a slightly more complex growth model, where the concentration of foam cells $\hat{c}_s = \hat{c}_s(\hat{x},t)$ ($\hat{x}\in{\hat{\SO}}$) is now governed by a PDE model, namely a non-stationary reaction-diffusion equation
\begin{equation} \label{PDEmodel}
\begin{aligned}
\partial_t \hat{c}_s - D_s \hat{\Delta} \hat{c}_s + R_s \hat{c}_s (1-\hat{c}_s)&= 0 \quad \text{on } {\hat{\SO}}, \; & D_s \hat{\partial}_n \hat{c}_s = \gamma(\sigmat_f^{WS}) &= \alpha \delta(x) \left(1+\frac{\left\|\sigmat_f^{WS}\right\|^2}{\sigma_0^{2}}\right)^{-1} \quad \text{on } \hat{\Gamma}, \\
\hat{c}_s &= 0 \quad \text{on } \hat{\Gamma}_s, \; & \hat{c}_s(\cdot,0) &= 0 \quad \text{on } \,{\hat{\SO}}, 
\end{aligned}
\end{equation}
where $\delta(x) = \min\left\lbrace 0, (x-1)(x+1)\right\rbrace^2$
and $\left\| \cdot \right\|$ is the Euclidean norm. Note that 
the norm $\|\sigmat_f^{WS}\|
|$ on the right-hand side of~\cref{PDEmodel} is a function in space and time, 
whereas the right-hand side of~\cref{reaction}, including $\|\sigmat_f^{WS}\|_{L^2(\Gamma)}$, depends only on time and is a constant in space. 
Furthermore, we have $\sigma_0$ and $ \alpha$ as in the ODE model~\cref{elasticgrowth-reference}, as well as positive diffusion and reaction coefficients $D_s$ and $R_s$, respectively.
The underlying idea is that the vessel wall is initially damaged in the central part of the interface around $\hat{x}_1=0$ and monocytes can penetrate into ${\hat{\SO}}$ in the part of the interface corresponding to $\hat{x}_1\in (-1,1)$, where $\delta$ is actually positive.

Moreover, in contrast to the first numerical example, we do not prescribe the shape of the plaque growth by a growth function $\hat g$. Instead, the growth part $F_g$ depends on the spatial distribution of the foam cell concentration via the relation
\begin{equation}\label{num:2:growth}
  \hat g(\hat x,t) =  1+ \hat{c}_s(\hat x,t),\quad 
  \hat\Ft_g(\hat x,t):=\hat g(\hat x,t)\,I.
\end{equation}
The Green Lagrange strain $\hat{\bs{E}}_e$ and the Piola--Kirchhoff stresses $\hat\Sigmat_e$ are then defined as above in \cref{Es,Piola}, with $\hat{g}$ given in~\cref{num:2:growth}. 

Note that, solving the PDE model has higher computational cost compared to the ODE model described in~\cref{sec:ODE}. On the other hand, a sophisticated reaction-diffusion problem may yield more realistic results for the plaque growth.
Nonetheless, such a model remains significantly cheaper compared to a fully coupled FSI problem; this discrepancy will only increase when moving to three-dimensional simulations with complex material models for the arterial wall. This observation is essential for the efficiency of our parallel time-stepping approach, as we will discuss in more detail in the following sections.

\section{Numerical framework}
\label{sec:num}

In this work, we use an Arbitrary Lagrangian Eulerian (ALE) approach to solve the FSI problem in~\cref{fullplaquemodel}. The ALE approach is the standard approach for FSI with small to moderate structural deformations;
see, e.g.,~\cite{DoneaSurvey,formaggia2010cardiovascular,RichterBuch,crosetto2011parallel,
bazilevs2013computational,wu2014fully,CauchaFreiRubio2018}. This assumption holds generally true for the simulation of plaque growth, unless the interest is to simulate a complete closure of the artery. To simulate a full closure, a Fully Eulerian formalism~\cite{Dunne2006, Richter2012b, FreiPhD} has been used in~\cite{FreiRichterWick2016, FreiPhD}; see also~\cite{AgerWalletal, BurmanFernandezFrei2020, BurmanFernandezFreiGerosa2022} for further works on FSI-contact problems in a Fully Eulerian or Lagrange-Eulerian formalism.

Given a suitable ALE map $\hat{\xi}_f =\hat{x} + \hat{\ut}_f$, its gradient $\hat{\Ft}_f = I + \hat{\nabla} \hat{\ut}_f$
and determinant $\hat{J}_f = {\rm det }\, \hat{\Ft}_f$, 
the ALE formulation of the FSI problem is given by:

\begin{varform}[Non-stationary FSI problem]\label{short-scale} 
  Find velocity $\hat \vt \in \vt^{\text{in}} + {\cal V}$,
  deformation $\hat\ut\in {\cal 
    W}$ and pressure $\hat p_f\in {\cal L}_f$, such that
  \[
  \begin{aligned}
    \big(\rho_f \hat J_f \hat \partial_t \hat
    \vt_f,\hat\phi\big)_{\hat\FL}
    + \big(\rho_f \hat J_f
    \hat\nabla\hat\vt_f\hat\Ft_f^{-1}(\hat\vt_f-\partial_t\hat\ut_f),
    \hat\phi\big)_{\hat\FL} 
    +
    \big(\hat\rho_s^0 \hat\partial_t \hat
    \ut_s,\hat\phi\big)_{\hat\SO}\qquad& \\
    + \big(\hat J_f\hat\sigmat_f\hat
    \Ft_f^{-T},\hat\nabla\hat\phi\big)_{\hat\FL}
    - \big( {\refone \hat J_f \nu_f \hat\Ft_f^{-1}\hat{\nabla}\hat{\vt}_f \hat\Ft_f^{-1} \hat{n}_f,
   \hat\phi}\big)_{\hat\Gamma^{\text{out}}}
    +\big(\hat\Ft_e(\hat{g})\hat\Sigmat_e(\hat{g}),\hat\phi\big)_{\hat\SO}
    &=0\;\forall \hat\phi\in {\cal V},\\
    \qquad\qquad\big(\widehat{\div}(\hat J_f\hat
    \Ft_f^{-1}\hat\vt_f),\hat\xi_f\big)_{\hat\FL}&=0\; 
    \forall \hat\xi_f\in {\cal L}_f,\\
    \big(d_t \hat \ut_s -\hat \vt_s,\hat\psi_s)_{\hat\SO} 
      &= 0 \; \forall \hat\psi_s \in {\cal L}_s,
  \end{aligned}
  \]
\end{varform}
{\refone 
\noindent where $\hat{\sigma}_f := \rho_f\nu_f (\hat\nabla \hat{\vt}_f \hat{\Ft}_f^{-1} + \hat{\Ft}_f^{-T} \nabla \hat{\vt}_f^T) -\hat{p}_f I$}.

There are different ways to compute the ALE map depending on the displacement $\hat{\ut}_s$ of the structural domain $\hat{\SO}$. In our simulations, we extend $\hat{\ut}_s$ into the fluid domain using harmonic extensions (called $\hat{\ut}_f$); that is, we solve a Laplacian problem with right hand side zero and $\hat{\ut}_s$ as Dirichlet boundary data on the fluid domain.

Note that $\hat\Sigmat_e$ and $\hat\Ft_e$ depend on $\hat{g}$ (and hence on $c_s$), as specified in \cref{Piola}. The function spaces are given by
\begin{align*}
  {\cal V} = [H^1_0(\hat\Omega;
    \hat\Gamma_f^{\text{in}}\cup\hat\Gamma_s)]^d,  \quad &{\cal W}   = [H^1_0(\hat\SO;
      \hat\Gamma_s)]^d, \quad
  {\cal L}_s:=L^2({\hat{\SO}}), \quad {\cal L}_f = L^2(\hat\FL).
\end{align*}
{\refone Furthermore, we assume that the solution $\hat\vt_f$ has higher regularity, such that the trace of $\hat\nabla \hat\vt_f$ is well-defined on $\hat{\Gamma}^{\text{out}}$, as needed in Variational formulation 1.}

\subsection{Temporal two-scale approach} \label{sec:two_scale}

Even for the simplified two-dimensional configuration considered in this work, 
a resolution of the micro-scale dynamics with a scale of milliseconds to seconds is unfeasible over the complete time interval of interest $[0,T_{\text{end}}]$, with $T_{\text{end}}$ being several months 
up to a year. For instance, when considering a relatively coarse micro-scale time step of $\delta \tau = 0.02$s, the number of time steps required to simulate a time frame of a whole year would be $365 \cdot 86\,400\cdot \frac{1s}{\delta\tau} \approx 1.58\cdot 10^9$, each step corresponding to the solution of a mechano-chemical FSI problem. 

This dilemma is frequently solved by considering a heuristic averaging: as the micro-scale is much smaller than the macro scale one considers a stationary limit of the FSI and solves for the stationary FSI problem on the macro scale (e.g., $\delta t\approx 1$ day); see, e.g.,~\cite{Chen2012, YangJaegerNeussRaduRichter2015,FreiRichterWick2016}. The wall-shear stress $\overline{\sigmat}_f^{WS}$ of the solution of this stationary FSI problem is then used to advance the foam cell concentration in~\cref{reaction} or~\cref{PDEmodel}.
  \begin{varform}[Stationary FSI problem]\label{long-scale} 
  Find velocity $\bar \vt \in \vt^{\text{in}} + {\cal V}$,
  deformation $\bar\ut\in {\cal 
    W}$ and pressure $\bar p_f\in {\cal L}_f$, such that
  \[
  \begin{aligned}
    \big(\rho_f \hat J_f
    \hat\nabla\bar\vt_f\bar\Ft_f^{-1}\bar\vt_f,
    \hat\phi\big)_{\hat\FL} 
    + \big(\bar J_f\bar\sigmat_f\bar
    \Ft_f^{-T},\hat\nabla\hat\phi\big)_{\hat\FL}\qquad&\\
    -\big( {\refone \hat J_f \nu_f \hat\Ft_f^{-1}\hat{\nabla}\hat{\vt}_f \hat\Ft_f^{-1} \hat{n}_f,
   \hat\phi}\big)_{\hat\Gamma^{\text{out}}}
    +\big(\bar\Ft_e(\hat{g})\bar\Sigmat_e(\hat{g}),\hat\phi\big)_{\hat\SO}
    &=0\;\;\forall \hat\phi\in {\cal W},\\
    \qquad\qquad\big(\widehat{\div}(\bar J_f\bar
    \Ft_f^{-1}\bar\vt_f),\hat\xi_f\big)_{\hat\FL}&=0\;\; 
    \forall \hat\xi_f\in {\cal L}_f
  \end{aligned}
  \]
\end{varform}
\noindent The foam cell concentration is then advanced by the ODE in~\cref{reaction} resp.\,the PDE in~\cref{PDEmodel}, with $\sigmat_f^{WS}$ replaced by $\overline{\sigmat}_f^{WS}$.
It has, however, been shown that $\gamma(\overline{\sigmat}_f^{WS})$ is not necessarily a good approximation of $\gamma(\sigmat_f^{WS})$, which depends on the pulsating blood flow; see the numerical results in~\cite{FreiRichterWick2016, Frei:2021:THN} and the analysis in~\cite{FreiRichter2020}. 
%
%

\begin{figure}
\begin{center}
\includegraphics[width=\textwidth]{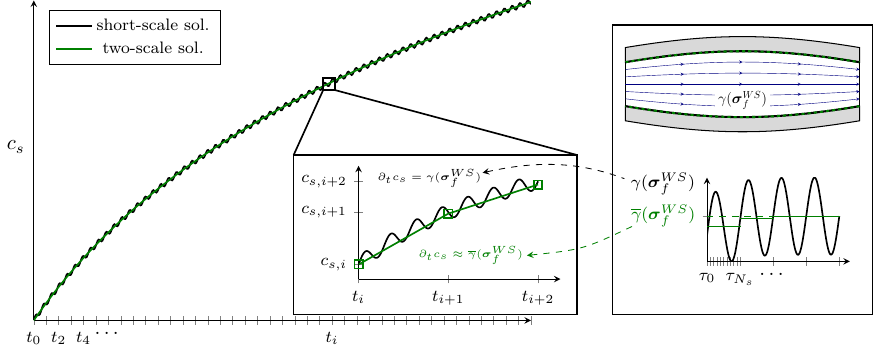}
\end{center}
\caption{Schematic representation of the two-scale algorithm as described in~\cref{sec:two_scale} and summarized in~\cref{twoscale}: instead of solving the FSI problem and updating the foam cell concentration $c_s$ on the micro time scale ($\delta\tau$, black), we compute an averaged growth function $\overline{\gamma}(\sigmat_f^{WS})$ over a few periods (heartbeats) of the micro-scale problem (right); this average is used to update $c_s$ on the macro scale ($\delta t$, green in the left part).
\label{fig:hom:parareal}}
\end{figure}

A more accurate two-scale approach has been presented by Frei and Richter in~\cite{FreiRichter2020}. The numerical approach can be cast in the framework of the Heterogeneous Multiscale Method (HMM); see, e.g.,~\cite{EEngquist2003,E2011,Abdulleetal2012}. 
In~\cite{FreiRichter2020}, a periodic-in-time micro-scale problem is 
solved in each time step of the macro scale, for instance each day. 
The growth function $\gamma(\sigmat_f^{WS})$ is then averaged by integrating over one period of the heart beat, its average will be denoted by $\overline{\gamma}(\sigmat_f^{WS})$. This averaged growth function is applied to advance the foam cell concentration by~\cref{reaction} or~\cref{PDEmodel}, respectively. A schematic illustration of the two-scale algorithm is given in~\cref{fig:hom:parareal}.

To be precise, we divide the macro-scale time interval $[0,T_{\text{end}}]$ into $N_f$ time steps of size $\delta t$
\begin{align}\label{deft}
0 = t_0 < t_1 < ... < t_{N_f} = {\sfrei T_{\text{end}}}, \qquad N_f = \frac{T_{\text{end}}}{\delta t}.
\end{align} 
As $c_s$ varies significantly on the macro scale only, using $c_s(t_m)$ as a fixed value for the growth variable on the micro scale results in a sufficiently good approximation in the time interval $[t_m,t_{m+1}]$. Then, one cycle of the pulsating blood flow problem (around $1$s) is to be resolved on the micro scale $\delta \tau$
\begin{align}\label{deftau}
0 = \tau_0 <\tau_1 <....< \tau_{N_s} = 1s, \qquad N_s = \frac{1s}{\delta \tau}.
\end{align}
It has been shown in~\cite{FreiRichter2020} (for a simplified flow configuration) that this approach leads to a model error ${\cal O}(\epsilon)$ compared to a full resolution of the micro scale, where $\epsilon=\frac{1s}{T_{\text{end}}}$ denotes the ratio between macro- and micro time scale and is in the range of ${\cal O}(10^{-7})$ for a typical cardiovascular plaque growth problem. The relations in~\cref{deft,deftau} imply $\delta \tau = \epsilon \delta t$.
In the model problems formulated above, the scale separation is induced by the parameter $\alpha={\cal O}(\epsilon)$.

A difficulty lies in solving the periodic micro-scale problem. If accurate initial conditions 
$\wt^0:=(\vt^0, \ut^0)$ are available on the micro-scale, a periodic solution can be computed using 
a time-stepping procedure for one cycle. If the initial conditions are known approximately, convergence to the periodic solution may still be obtained after simulating a few cycles of the micro-scale problem due to the dissipation of the flow problem; see~\cite{FreiRichter2020,Frei:2021:THN}. Numerically, it can be checked after each cycle if the solution is sufficiently close to a periodic state. In this work, we apply a stopping criterion based on the computed averaged growth value:
\begin{align*}
|\overline{\gamma}(\sigmat_f^{WS,r}) - \overline{\gamma}(\sigmat_f^{WS,r-1})| < \epsilon_p,
\end{align*}
where $r=1,2,...$ denotes the iteration index with respect to the number of cycles of the micro-scale problem. The algorithm is summarized as~\cref{twoscale}, where we use the abbreviation $\wt^{r,s}:= (\vt^{r,s}, \ut^{r,s})$.



\begin{algorithm}
\caption{Two-Scale Algorithm\label{twoscale}}
  Set suitable starting values $\wt^{0,0}=(\vt^{0,0}, \ut^{0,0})$ and time-step sizes $\delta t, \delta T$.
  \\
  \For{$n=1,2,\dots, N_f$}{ 
  \begin{enumerate}
    \item[1.)] \textbf{Micro problem}: Set $r\gets 0$\\
    \While{$|\overline{\gamma}(\sigmat_f^{WS,r}) - \overline{\gamma}(\sigmat_f^{WS,r-1})| > \epsilon_p$}{
    \begin{itemize}
  \item[1.a)]
  Solve micro-scale problem in~\cref{fullplaquemodel} 
    in $I_n =(t_n,t_n+ 1s)$
    \[    
    \{\wt^{r,0},c_s^{n-1}\}\mapsto 
    \{\wt^{r,m}\}_{m=1}^{N_s}
    \]
  \item[1.b)] Compute the averaged growth function
    \[
    \overline{\gamma}(\sigmat_f^{WS,r+1}) = \frac{1}{N_s}\sum_{m=1}^{N_s}
    \gamma(\sigmat_f^{WS,m}(\vt^{r+1,m}), c_s^{n-1})
    \]
    and set $\wt^{r+1,0}= \wt^{r,N_s}, \, r\gets r+1$.
    \end{itemize}}
  \item[2.)] \textbf{Macro problem}: Update the foam cell concentration $c_s^n$ by~\cref{reaction} or~\cref{PDEmodel}.
    \end{enumerate}\vspace{-0.4cm}
}
\end{algorithm}

As starting values $\wt^{r,0}=(\vt^{r,0}, \ut^{r,0})$ in step 1.~of the algorithm, we use the variables $\wt^{r-1,N_s}=(\vt^{r-1,N_s},\ut^{r-1,N_s})$ from the quasi-periodic state of the previous macro step. It has been observed in~\cite{Frei:2021:THN} that these are usually closer to the starting values of the periodic state than the solution of an averaged stationary problem on the macro scale.

\subsection{Numerical example}
\label{sec:numex1}

Before we present the different approaches for parallel time-stepping, let us illustrate the two-scale approach by a first numerical example. Therefore, we use the simple ODE growth model introduced in~\cref{sec:ODE},~\cref{reaction}. The test configuration introduced here will be used in~\cref{sec:parallel} as well.

Concerning the geometry we use a two-dimensional channel of length $10$\,cm and an initially constant width $\omega(0)$ of $2$\,cm as illustrated in~\cref{fig:conf}. 
The solid parts on the top and bottom corresponding to the arterial wall have an initial thickness of $1$\,cm each.
Fluid density and viscosity are given by $\rho_f=1\,\text{g}/\text{cm}^3$ and
$\nu_f=0.04\,\text{cm}^2/\text{s}$, respectively. The growth parameters are set to $\sigma_0=30 \frac{\text{g\,cm}}{\text{s}^2}$ and $\alpha=5\cdot 10^{-7}\frac{1}{\text{s}}$, which yields a realistic time-scale for the arterial plaque growth, see~\cite{FreiRichter2020}.
The solid parameters are
$\rho_s=1\,\text{g}/\text{cm}^3$, $\mu_s=10^4\,{\refone \text{dyne}/\text{cm}^2}$, and
$\lambda_s=4\cdot 10^4\,\text{dyne}/\text{cm}^2$. As an inflow boundary condition, we prescribe a pulsating
velocity inflow profile on $\Gamma_f^\text{in}$ 
given by
\begin{equation}\label{num:1:inflow}
  \vt_f^\text{in}(t,x) = 30\begin{pmatrix}
   \sin(\pi {\refone \frac{t}{\cal P}})^2 (1-x_2^2)\\ 0
  \end{pmatrix}  \text{cm}/\text{s}.
\end{equation}
{\refone Here, ${\cal P} = 1$\,s is the period of a single heartbeat.}
The symmetry of the configuration can be exploited in order to reduce computational cost by simulating only on 
the lower half of the computational domain 
and imposing the symmetry condition $\vt_f\cdot \nt=0,\, \tau^T\sigmat_f \nt =0$ on the symmetry boundary $\Gamma_{\text{sym}}$, where $\tau$ denotes a tangential vector.

We discretize the FSI problem (\cref{short-scale}) in time using the backward Euler method. For discretizing the ODE growth model, we use the forward Euler method, 
which results in 
    \begin{align}\label{eq:EulerODE}
    c_s^{n} = c_s^{n-1} + \delta t\, \overline{\gamma}(\sigmat_f^{WS}, c_s^{n-1}), \qquad n=1,...,N_f.
    \end{align}

For spatial discretization, we use biquadratic ($Q_2$) equal-order finite elements for all variables and LPS stabilization~\cite{BeckerBraack2001} for the fluid problem. 
Our mesh, containing both fluid and solid, consists of 160 rectangular grid cells; this corresponds to a total of 3\,157 degrees of freedom. The time-step sizes are chosen as $\delta \tau=0.02$\,s and $\delta t=0.3$ days; the tolerance for periodicity of the micro-scale problem as $\epsilon_p=10^{-3}$. All the computational results have been obtained with the finite element library Gascoigne3d~\cite{Gascoigne}. We use a fully monolithic approach for the FSI problem following Frei, Richter \& Wick~\cite{FreiRichterWick2016, RichterBuch}.

\begin{figure}[t]
\center
\includegraphics[width=0.49\textwidth]{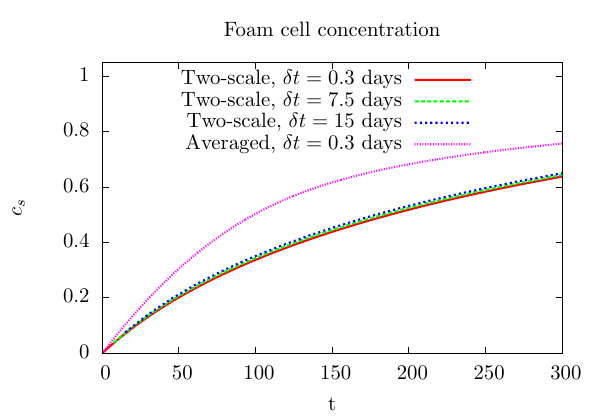}
\includegraphics[width=0.49\textwidth]{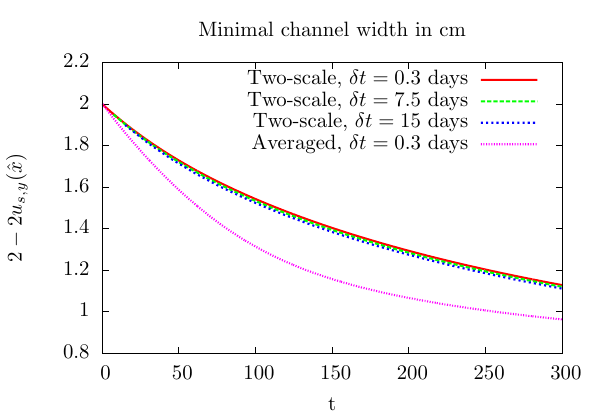}\\
\includegraphics[width=0.5\textwidth]{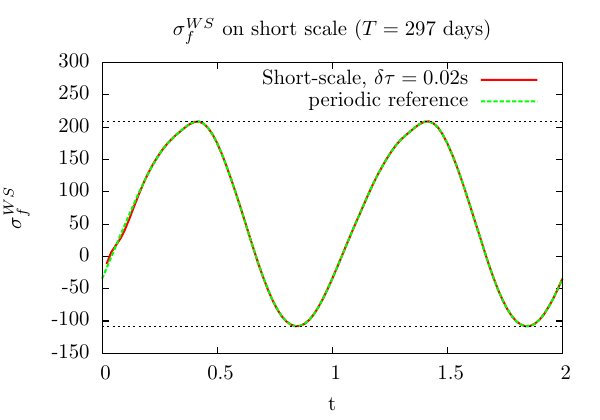}
\caption{\label{fig:NumEx} \textbf{Top:} Comparison of a pure averaging approach with the two-scale approaches with different macro-scale time steps $\delta T$. \textbf{Top left}: Concentration $c_s$ over time. \textbf{Top right}: Channel width over time. \textbf{Bottom}: Wall shear stress $\sigmat_f^{WS}$ over two periods of the micro scale at $t=297$ days for $\delta t=9$ days. The initial values are taken from the previous micro problem at $t=288$ days. 
}
\end{figure}

In~\cref{fig:NumEx}, we compare results obtained with the two-scale approach for different macro time-step sizes $\delta t$ with the heuristic averaging approach outlined above; see~\cref{long-scale}. We see that the pure averaging approach underestimates the growth significantly. Even a very coarse discretization of the macro-scale time interval $\delta t=15$ days in the two-scale approach gives a much better approximation. We note, however, that 
too coarse time steps
might introduce different issues, as the starting values $\vt^{0,0}$ and $\ut^{0,0}$ might not be good approximations for the periodic state on the micro scale anymore. In combination with a significant mesh deformation for $t\approx 300$ days, this led to divergence of~\cref{twoscale} for an even coarser macro time-step size $\delta t=30$ days. For $\delta t \leq 15$ days, a near-periodic state was reached in 2-3 iterations in each macro-step. A detailed convergence study for a very similar problem has been given in~\cite{Frei:2021:THN}.

At the bottom of~\cref{fig:NumEx}, we show the $L^2(\Gamma)$-norm of the wall shear stress 
over two periods, i.e., 2 seconds, of the micro problem for $\delta t=0.02$s and $\delta T=9$ days. The initial values are taken from the periodic state of the micro problem at the previous macro time-step, i.e., 9 days before. We see that the wall shear stresses converge very quickly to the periodic state. An initial deviation of $6.28 \frac{\text{g cm}}{\text{s}^2}$ to the periodic reference solution at time $\tau=0.02$ is reduced to $6\cdot 10^{-4} \frac{\text{g cm}}{\text{s}^2}$ at time $\tau=1$ within only one period.

\section{Parallel time-stepping}
\label{sec:parallel}

The main cost in~\cref{twoscale} lies in the solution of the non-stationary micro problem in step 1.a, which needs to be solved in each time step of the macro problem. Considering a relatively coarse micro-scale discretization of $\delta t = 0.02$\,s, as used in the previous section, 50 time steps are necessary to compute a single period of the heart beat. The simulation of two or more cycles might be necessary to obtain a near-periodic state in step 1.a of~\cref{twoscale}; cf. the discussions in~\cite{FreiRichter2020,Frei:2021:THN}. In a realistic scenario, each time step of the micro problem corresponds to the solution of a complex three-dimensional FSI problem, which makes already the solution of one micro problem very costly.

For this reason, parallelization needs to be exploited in different ways: First, one can make use of spatial parallelization using a scalable solver, for instance, based on domain decomposition or multigrid methods; see, e.g.,~\cite{Heinlein:2016:PIT,gee2011truly,kong2018scalability,wu2014fully,barker2010scalable,deparis2016facsi,heinlein2016parallel}. 
Since the focus of this contribution is on parallelization in time, we will not discuss this aspect here. In particular, we expect that even when the speedup due to spatial parallelization saturates, the computing times for the whole plaque growth simulations will remain unfeasible. Therefore, we have to make use of an additional level of parallelization, that is, temporal parallelization of the macro-scale problem. In particular, we will use the parareal algorithm, which we will recap in the next subsection.

In order to motivate the algorithmical developments in this section, we will already 
present some first numerical results for the ODE growth model in~\cref{reaction} within the section.

\subsection{The parareal algorithm} \label{sec:parareal}

First, the time interval of interest $[0,T_{\text{end}}]$ on the macro scale is divided into $P$ sub-intervals $I_p=[T_{p-1}, T_p]$ of equal size, where
\begin{align}\label{eq:coarse_disc}
0 = T_0 < T_1 < ... < T_P = T_{\text{end}}.
\end{align}

In order to define the parareal algorithm, suitable fine and coarse problems need to be introduced. Note again that we will apply the parareal algorithm only on the macro scale as introduced in~\cref{sec:two_scale}; hence, both the fine and the coarse scale of parareal correspond to the macro scale of the homogenization approach. The fine problem advances the growth variable $c_s$ from time $T_p$ to $T_{p+1}$ by solving~\cref{twoscale} with a smaller time-step size $\delta t$ (e.g., $0.3$ days) on the corresponding fine
time discretization of $[T_p, T_{p+1}]$:
\begin{align*}
T_p = t_{p,0} < t_{p,1} < ... < t_{p,n_p} = T_{p+1}, \qquad n_p = \frac{T_{p+1}-T_p}{\delta t}, \qquad  t_{p,q} := t_{p\cdot n_p +q} = t_{p,q-1} + \delta t.
\end{align*}
The fine propagator on a process $p$ consists of a time-stepping procedure to advance $c_s(T_p)$ to $c_s(T_{p+1})$ with the fine time step $\delta t$. We write
\begin{align*}
{\sfrei c_s^{\text{fine}}}(T_{p+1}) = {\cal F}({\sfrei c_s^{\text{fine}}}(T_p)).
\end{align*}

The efficiency of the parareal algorithm depends strongly on the computational cost of the coarse propagator. In particular, it needs
to be much cheaper than the fine problems since it is defined globally on $[0,T_{\text{end}}]$ and, hence, introduces synchronization. {\reftwo Thus, we use a large time-step $\delta T$ and
\begin{align*}
T_p = \bar{T}_{p,0} < ... < \bar{T}_{p,N_p} = T_{p+1}, \qquad N_p = \frac{T_{p+1}-T_p}{\delta T}, \qquad  \bar{T}_{p,q} := \bar{T}_{p\cdot N_p +q} = \bar{T}_{p,q-1} + \delta T.
\end{align*}
For simplicity, we assume that both time-step sizes $\delta t$ and $\delta T$ are uniform throughout $[0,T_\text{end}]$.
In order to keep the cost of the coarse propagator as low as possible, we will mainly focus on the case that the coarse time steps coincide with the $P$ sub-intervals $I_p$ in the numerical results, i.e., $N_p=1$, such that the total number of coarse time steps
\begin{align*}
N_c:= P\cdot N_p
\end{align*}
is equal to $P$.}
We denote {\reftwo the coarse propagation from $T_p$ to $T_{p+1}$} by
\begin{align*}
{\sfrei c_s^{\text{coarse}}}(T_{p+1}) = {\cal C}({\sfrei c_s^{\text{coarse}}}(T_p)).
\end{align*}

We use capital letters $T_p$ to denote the coarse discretization of $[0,T_{\text{end}}]$ into $P$ parts {\reftwo and for the time-steps $T_{p,q}$ on the coarse level.
By} small letters $t_{p,q}$ we denote the finer discretization on {\reftwo each sub-interval;} 
the two discretizations yield the first level (fine problem) and second level (coarse problem) of the parareal algorithm. {\reftwo Of course, it is also possible to employ coarse time steps which differ from the sub-intervals on the fine level, but for the sake of simplicity, we will not discuss this case in this work.}
In the example given above with $T_{\text{end}}=300$ days, $\delta T$ would be 30 days for {\reftwo $P=N_c=10$}, while $\delta t$ is 0.3 days. On the micro scale, the times $\tau_i$ and time-step size $\delta\tau$ are defined locally in $[t_i, t_i + 1s]$; in the example above, we had $\delta \tau =0.02$s.
Note that the micro scale influences the parareal algorithm only indirectly due to the temporal homogenization approach.

\begin{figure}
\begin{center}
\includegraphics[width=0.55\textwidth]{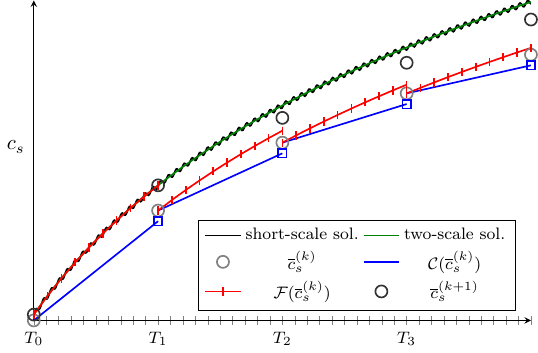}
\end{center}
\caption{Schematic representation of the parareal algorithm applied to the homogenized (two-scale) problem as shown in~\cref{fig:hom:parareal}: micro-scale solution (black), two-scale solution (green), parareal iterates in iteration $k$ (light grey) and $k+1$ (dark grey), fine scale solution (red), and coarse scale solution (blue).
\label{fig:parareal}}
\end{figure}

\begin{algorithm}[t]
	\caption{Parareal algorithm\label{alg:parareal}}
	\begin{enumerate}
		\item[(I)] \textbf{Initialization}: Compute $\big\{(\overline{c}_s^{(0)}(T_p), \wt^{(0)}(T_p)\big\}_{p=1}^{P}$ by means of~\cref{twoscale} 
		with a coarse macro time-step size $\delta T := {\reftwo (T_{p}- T_{p-1}) / N_p}$. 
		Set $k\gets 0$\\
		\item[(II)]\While{$|c_s^{(k+1), {\sfrei \text{fine}}}({\sfrei T_{\text{end}}}) - c_s^{(k), {\sfrei \text{fine}}}({\sfrei T_{\text{end}}})| > \epsilon_{\text{par}}$}{
			\begin{itemize}
				\item [(II.a)] \textbf{Fine problem}: \\
				\For{$p=1,...,P$}{ 
					\begin{itemize}
						\item[(i)] Initialize $c_s^{(k+1), {\sfrei \text{fine}}}(T_p)=\overline{c}_s^{(k)}(T_p), \,\wt^{(k+1)}(T_p)  = \wt^{(k)}(T_p)$
						\item[(ii)] Compute $\big\{(c_s^{(k+1), {\sfrei \text{fine}}}(t_{p,q}), \wt^{(k+1)}(t_{p,q})\big\}_{q=1}^{n_p}$ by~\cref{twoscale} with fine time-step size\\ $\delta t$ and set 
						${\cal F}(\overline{c}_s^{(k)}(T_p)) = c_s^{(k+1), {\sfrei \text{fine}}}(t_{p, n_p})$
					\end{itemize}  
				}
				\item[(II.b)] \textbf{Coarse problem}\\ 
				\For{$p=1,...,P$}{ 
					\begin{itemize}
						\item[(i)] Compute ${\cal C} (\overline{c}_s^{(k+1)}(T_p))$ by solving one time step of~\cref{twoscale} with time-step\\ size $\delta T={\reftwo (T_{p}- T_{p-1}) / N_p}$
						\item[(ii)] Parareal update
						\begin{align*}
							\overline{c}_s^{(k+1)}(T_{p+1}) = {\cal C} (\overline{c}_s^{(k+1)}(T_p)) +{\cal F}( \overline{c}_s^{(k)}(T_p)) - {\cal C}(\overline{c}_s^{(k)}(T_p)).
						\end{align*}
				\end{itemize}}\vspace{-0.3cm}
			\end{itemize}
			$k\gets k+1$}
	\end{enumerate}
\end{algorithm}

Then, given an iterate $\overline{c}_s^{(k)}$ for some $k\leq 0$, the parareal algorithm computes $\overline{c}_s^{k+1}$ by setting
\begin{align}\label{eq:parareal}
\overline{c}_s^{(k+1)}(T_{p+1}) = {\cal C} (\overline{c}_s^{(k+1)}(T_p)) +{\cal F}( \overline{c}_s^{(k)}(T_p)) - {\cal C}(\overline{c}_s^{(k)}(T_p)) \quad \text{ for }\, p=0,...,P-1.
\end{align}
This can be seen as a predictor-corrector scheme, 
where the coarse predictor ${\cal C} (\overline{c}_s^{(k+1)}(T_p))$ is corrected by fine-scale contributions that depend only on the previous iterate $\overline{c}_s^{(k)}$; thus, the fine problems {\sfrei can} be computed fully in parallel.
A schematic illustration of the parareal algorithm is given in~\cref{fig:parareal}.

Let us analyze the application of~\cref{eq:parareal} to the two-scale problem (\cref{twoscale}) in more detail. The first term in~\cref{eq:parareal} requires the solution of one micro problem and an update of the foam cell concentration (by~\cref{reaction} or~\cref{PDEmodel}) in each coarse time step $T_p\to T_{p+1}$. 
Within the fine-scale propagator (second term in~\cref{eq:parareal}) $n_p=\lceil N_f/P\rceil$ time-steps need to be computed per process, where $\lceil g \rceil$ denotes the next-biggest natural number to $g$ ($\widehat{=}$ \texttt{ceil(g)}) {\refone and $N_f$ is the number of macro time steps; see~\cref{deft}}.  Each time step requires the solution of one micro problem and an update of the foam cell concentration. The last term in~\cref{eq:parareal} has already been computed in the previous iteration (compare the first term on the right-hand side of the same equation) and thus introduces no additional computational cost. The algorithm is summarized in~\cref{alg:parareal}.

%
{\reftwo We use the FSI variables $w^{(k)}(T_p)$ from the previous parareal iteration 
as initial values in step (II.a)(i). 
The initial values $\overline{c}_s^{(k)}(T_p)$ in step (II.a)(i) 
are taken from step (II.b) of the coarse problem.
Moreover, to initialize the variables $\{\overline{c}_s^{(0)}(T_p), \wt_s^{(0)}(T_p)\}_{p=1}^P$ before the first parareal iteration, we apply the coarse propagator once
with a large time-step size $\delta T= (T_p-T_{p-1}) / N_p$.}

For the ODE growth model, we use the value of ${\sfrei c_s^{\text{fine}}}$ at the end time $T_{\text{end}}$
\begin{align}\label{eq:terminal}
|c_s^{(k+1),{\sfrei \text{fine}}}(T_{\text{end}}) - c_s^{(k), {\sfrei \text{fine}}}(T_{\text{end}})| \leq \epsilon_{\text{par}}. 
\end{align}
as the stopping criterion for the parareal algorithm.

{\reftwo
\subsubsection{Parallelization approaches} 
\label{sec:parallelization}

\begin{figure}[t]
	\begin{center} 
		Master-slave parallelization \\
		\includegraphics[width=\textwidth]{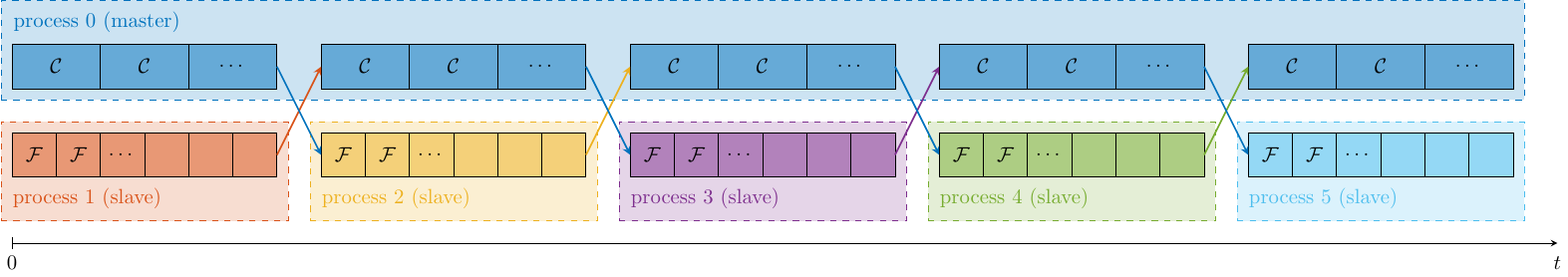} \\[2mm]
		Distributed parallelization \\
		\includegraphics[width=\textwidth]{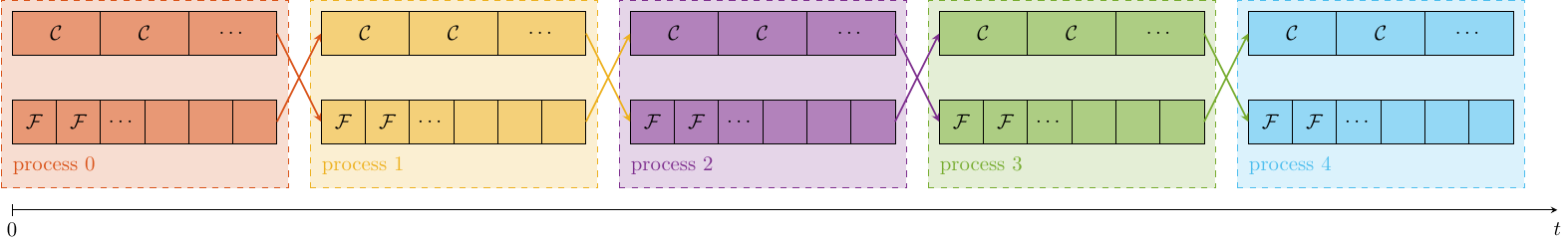}
	\end{center}
	\caption{\reftwo Master-slave and distributed parallelization schemes for the parareal approach. In the master-slave case, the coarse propagation is performed on a separate process (blue) and information is communicated between master and slave. In the distributed case, the coarse propagation and the fine propagation is performed on the same processes. Communication is indicated by an arrow in color of the sender.
	\label{fig:parallelization}
	}
\end{figure}

As mentioned earlier, in this work, we focus on studying the effectivity and efficiency of our approach by investigating the convergence of the parareal algorithm based on a serial implementation. Since the computational cost of the micro FSI problems is very high, we expect the communication cost to be negligible; cf.~\cref{sec:compcost} as well as~\cref{sec:theo_cost,sec:runtimes} for a more detailed discussion of the computational cost for the ODE and the PDE growth model, respectively. Even though our implementation used in the numerical examples below is only serial, we want to briefly discuss two potential parallelization approaches for a parallel implementation of~\cref{alg:parareal}: \textit{master-slave parallelization} and \textit{distributed parallelization}; cf.~\cref{fig:parallelization} for a sketch of both approaches for our parareal algorithm. Note that, as also indicated by~\cref{fig:parallelization}, we always assume a one-to-one correspondence of fine problems and processes in this paper.

In the master-slave parallelization approach, the coarse problem is assigned to a single process and the local problems on the fine level are distributed among the remaining processes. 
Each slave process $p$ has to communicate ${\cal F}(\overline{c}_s^{(k)}(T_p))$ to the master after computation of the fine problem for the parareal update (step II.a in~\cref{alg:parareal}), and the master process has to communicate $\overline{c}_s^{(k+1)}(T_p)$ back after the parareal update (step II.b in~\cref{alg:parareal}); this corresponds to all-to-one and one-to-all communication patterns.

In the distributed parallelization approach, the parallelization of the coarse problem is different. In particular, time intervals are assigned to the different processes, and each process computes both the fine problem and the part of the coarse problem on those time intervals; for simplicity, we do not discuss the case that the partition of the coarse problem onto the processes is different from the local (fine) problems, which would also be possible. As can be seen in~\cref{fig:parallelization}, a major benefit of this approach is the different communication pattern: instead of all-to-one and one-to-all communication, each process sends the parareal update $\overline{c}_s^{(k+1)}(T_p)$ to the next process in line. Note that, in terms of the coarse problem, the processes can only perform their computations in serial, whereas the local problems can, again, be computed concurrently.

In some implementations of two-level methods, the coarse problem is assigned to one of the processes dealing with the local problems instead of an additional process; if the memory and computational cost of the coarse problem is negligible, this can be beneficial since all available processes can be employed for the local problems. It can be seen as a mixture of the master-slave and distributed parallelization approaches. In our case, this might not be feasible since the coarse problem is defined on the same spatial mesh as the local problems, leading to significant memory cost. Furthermore, more sophisticated approaches such as task-based scheduling~\cite{aubanel_scheduling_2011} can be employed for the parallelization.

An investigation of the differences in performance of the different approaches can only be performed using an actual parallel implementation. This is out of the scope of this paper, and therefore, we will leave it to future work. 
}

\subsubsection{Computational costs}
\label{sec:compcost}

Even in the simplified two-dimensional configuration considered in this work,
it suffices to count the number of micro problems to be solved to estimate the computational cost of corresponding parallel computations and to
compare 
different algorithms. Remember that, in each micro problem, $r_p\cdot \frac{1}{\delta t}$ FSI systems need to be solved, where $r_p$ is the number of cycles required to reach a near-periodic state. In the example considered above, these were at least $2\cdot \frac{1}{0.02}=100$ costly FSI problems.  
Hence, our discussion in the following sections will be based on the assumption that the computational cost for all other  steps of the parareal algorithm as well as the communication between the processes can be neglected. This is particularly obvious for the ODE growth model in~\cref{reaction}. For the PDE growth model in~\cref{PDEmodel}, we will discuss the computational and communication cost in more detail 
in~\cref{sec:theo_cost,sec:runtimes}. 

To analyze the computational cost of~\cref{alg:parareal} under this assumption, let us denote the total number of iterations of the
parareal algorithm by $k_{\text{par}}$. {\reftwo As discussed in~\cref{sec:parareal,sec:parallelization}, we assume that the number of coarse time steps is the same as the number of fine problems and processes, respectively. Therefore,}
we need to solve {\reftwo $N_c$} micro problems in step I, $k_{\text{par}}\cdot \lceil N_f/P\rceil$ micro problems on each of the $P$ processes in step II.a (ii), and $k_{\text{par}}\cdot {\reftwo N_c}$ micro problems {\reftwo within the coarse propagation} in step II.b (i). This corresponds to the solution of 
\begin{align}\label{eq:comp_cost}
	\underbrace{k_{\text{par}}\cdot \lceil N_f/P \rceil}_{\text{fine level ($P$ parallel processes)}} + \underbrace{(k_{\text{par}}+1) \cdot {\reftwo N_c}}_{\text{coarse level (1 serial process)}}
\end{align}
serial micro problems{\reftwo; this is the case for both parallelization approaches discussed in~\cref{sec:parallelization}.} {\reftwo If the coarse time step $\delta T$ is chosen independently of $T$, $N_c$ is fixed and the computational cost tends to saturate for large $P$
(at least if we assume that the number of required parareal iterations $k_{\text{par}}$ is independent of $P$). The cheapest possible coarse propagator, on the other hand, is to use one coarse time step per coarse propagation, i.e., $N_p=1$ and $N_c=P$. In this case the computational time increases for large $P$, if we assume that the number of parareal iterations $k_{\text{par}}$ is independent of $P$. The minimum computational time is attained for $P\approx \sqrt{N_f}$.}

{\reftwo As a serial implementation requires the solution of $N_f$ micro problems, the speed-up 
of the proposed parareal algorithm is given by
\begin{align}\label{eq:speedup}
\text{speedup}(P) = \frac{N_f}{k_{\text{par}}\cdot \lceil N_f/P \rceil  + (k_{\text{par}}+1) \cdot {N_c}} \approx \frac{1}{\frac{k_{\text{par}}}{P}  + (k_{\text{par}}+1) \frac{N_c}{N_f}},
\end{align}
if we assume a perfect load balancing on the fine processes. This is a special case of the standard analysis for the speed-up of the parareal algorithm; see, e.g.,~\cite{Ruprecht2017}.
}

The aim of this work is to show potential for the parareal algorithm as a parallel time-stepping method for plaque growth simulations. A final assessment of the parallelization capabilities can, of course, only be made based on computing times of an actual parallel implementation. {\reftwo We will leave this to future work.}


{\reftwo \subsubsection{Convergence theory}
\label{subsec.conv}

Let us briefly recap the standard convergence theory for the parareal algorithm; see, e.g.,~\cite{GanderHairer}. For this purpose, let us consider the simpler ODE model
\begin{equation}\label{eq:EulerODEcont}
    \partial_t c_s = \gamma(c_s), \qquad c_s(0)=0.
\end{equation}
and its discretization by the explicit Euler method
\begin{align}\label{eq:EulerODE2}
c_s^{n} = c_s^{n-1} + \delta T\, \gamma(c_s^{n-1}), \qquad c_s^0=0.
\end{align}
We assume that the function $\gamma$ and its derivative are Lipschitz continuous in $c_s$, i.e., we assume 
\begin{align}\label{Lipschitz}
|\gamma(c_s^1) - \gamma(c_s^2)| \leq L |c_s^1 - c_s^2|, \quad
\big|\frac{d}{dt} \gamma(c_s^1) - \frac{d}{dt} \gamma(c_s^2) \big| \leq L |c_s^1 - c_s^2| \quad \text{for all } c_s^1, c_s^2 \in\mathbb{R}
\end{align} 
for some $L>0$.
As in \cite{GanderHairer}, we also assume for simplicity that the fine-scale propagator advances the ODE exactly, i.e.,
\begin{align}\label{exact}
c_s(T_n)  = {\cal F}(c_s(T_{n-1})) = c_s(T_{n-1}) + \int_{T_{n-1}}^{T_n} \gamma(c_s(s))\, \text{d}s. 
\end{align}
This is motivated by the fact that the time discretization error in the fine propagator is typically small compared to the coarse propagator.
The following convergence result is shown by Gander \& Hairer in~\cite{GanderHairer}:

\begin{lemma}\label{lem.conv_parareal}
Let $c_s^{*}\in C^1(0,T_\text{end})$ 
be the exact solution of~\cref{eq:EulerODEcont}, and let $\{\overline{c}_s^{(k)}(T_n)\}_{n=1}^{P}$ be the $k$-th iterate of the parareal algorithm in~\cref{eq:parareal} using the forward Euler method in~\cref{eq:EulerODE2} as the time discretization. Under the assumptions made above, it holds for the error $ e_n^{(k)} = |\overline{c}_s^{(k)}(T_n)-c_s^{*}(T_n)|$ and $k\in\mathbb{N}_0, n\in \{1,\ldots,P\}$ that
\begin{align}
e_n^{(k)} 
&\leq   \left(cT_n\right)^{k+1} \beta^{n-k-1} \delta T^{k+1} \max_{t\in [0,T_\text{end}]} |\partial_t c_s^{*}(t)|\label{eq:est_parareal}.
\end{align}
with a constant $c>0$
and $\beta=1+L\delta T$.

\end{lemma}

\begin{remark}
    The assumptions in~\cref{Lipschitz} cannot be easily verified for the ODE model in~\cref{reaction}, as the fluid forces and $\bar{\gamma}(\sigma_f^{WS})$ depend in a highly nonlinear way on $c_s$. In the numerical examples given below, we will, however, observe that the convergence behavior is similar to the one shown in~\cref{lem.conv_parareal}.
\end{remark}

}

\subsubsection{Numerical results}
\label{sec:pararealres}

\begin{figure}[t!]
\centering
\includegraphics[width=0.48\textwidth]{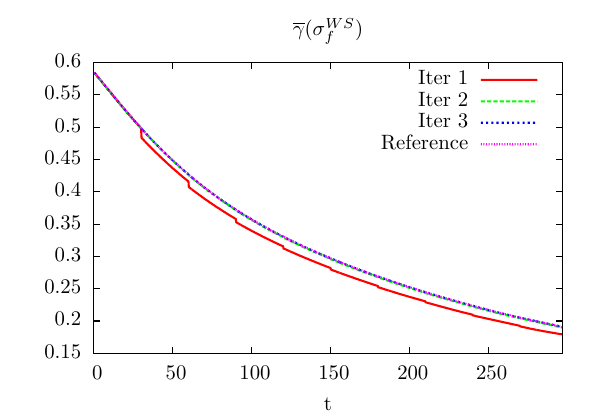}
\includegraphics[width=0.48\textwidth]{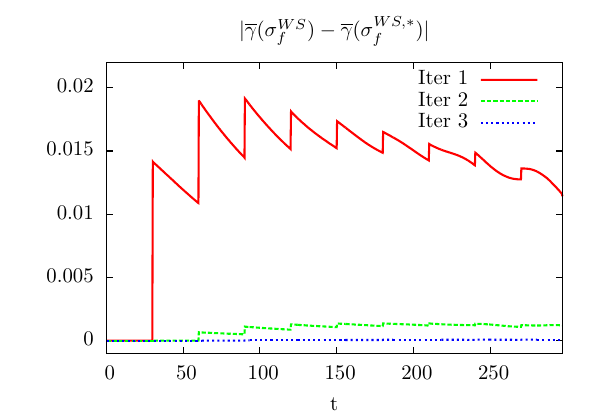}\\
\includegraphics[width=0.48\textwidth]{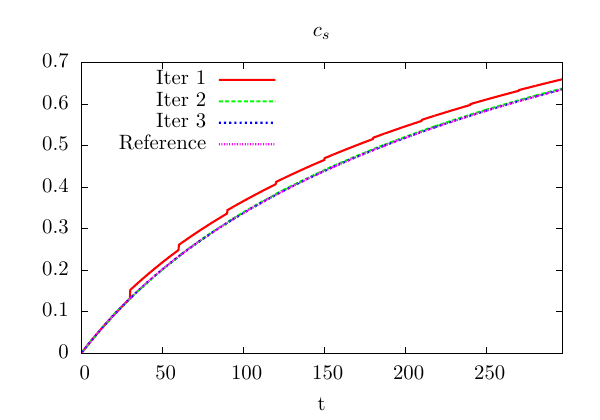}
\includegraphics[width=0.48\textwidth]{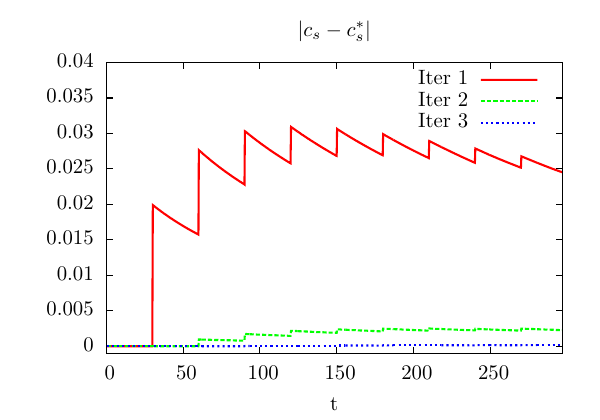}\\
\includegraphics[width=0.48\textwidth]{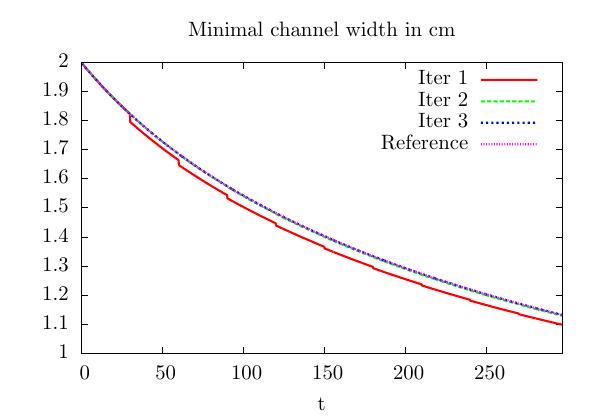}
\includegraphics[width=0.48\textwidth]{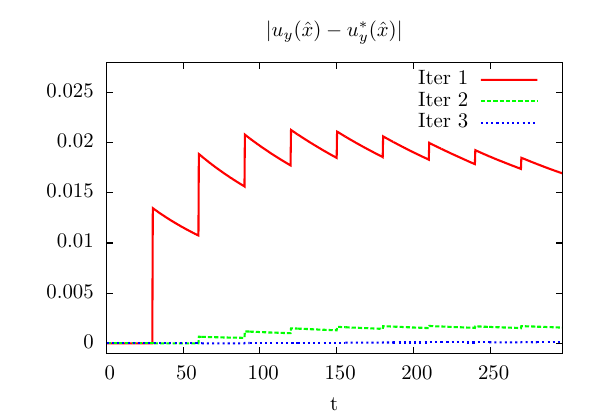}\\
\caption{\label{fig:parallel} Behavior of the first three iterates of the parareal algorithm in the first numerical example (ODE growth) for $P=10$. A description of the problem is given in~\cref{sec:pararealres}, numerical values in~\cref{tab:csk}. The iterates converge very quickly towards the reference solution. \textbf{Top}: Growth function {\reftwo and error w.r.t the reference solution $\sigmat_f^{WS}(\vt_f^*)$} over time $t$. \textbf{Center}: {\sfrei $c_s^{(k),\text{fine}}$} {\reftwo and error w.r.t the reference solution $c_s^*$} over time. \textbf{Bottom}: Channel width over time {\reftwo and error $|\ut_{s,y}(\hat{x}) - \ut_{s,y}^*(\hat{x})|$ w.r.t. the vertical component of the reference solution $\ut_s^*$ at the narrowest point $\hat{x}$} over time .}
\end{figure}

We use, again, the simple example described in~\cref{sec:numex1} with the ODE growth model in~\cref{reaction} and set $\epsilon_{\text{par}}=10^{-3}$.
Results for $P=N_c=10$ and $\delta t=0.3$ days are shown in~\cref{fig:parallel}, where the first 3 iterates are compared against a reference solution {\reftwo $(c_s^*, \ut_s^*, \vt_f^*)$, }
which was computed by a standard serial time-stepping scheme,
as in~\cref{sec:numex1}.
We observe fast convergence towards the reference curve {\reftwo in all three quantities of interest}. The stopping criterion in~\cref{eq:terminal} was satisfied after 4 iterations of the parareal algorithm.

\begin{table}[t]
\begin{center} \small
\begin{tabular}{|l|rrrrr|r|}
\hline
k &$P=10$ &$P=20$ &$P=30$ &$P=40$ &$P=50$ & ref. (serial) \\
\hline
1 &$2.21\cdot 10^{-2}$ &$1.19\cdot 10^{-2}$ &$8.12\cdot 10^{-3}$
&$5.83\cdot 10^{-3}$ &$5.46\cdot 10^{-3}$ &{\reftwo 0}\\
2 &$2.24\cdot 10^{-3}$ &$5.63\cdot 10^{-4}$ &$2.62\cdot 10^{-4}$ &$1.37\cdot 10^{-4}$ &$8.61\cdot 10^{-5}$ &-\\
3 &$1.42\cdot 10^{-4}$ &$2.02\cdot 10^{-5}$&$6.53\cdot 10^{-6}$ &$2.32\cdot 10^{-6}$ &$8.30\cdot 10^{-7}$ &-\\
4 &$5.76\cdot 10^{-6}$ &- &- &- &- &-\\
\hline 
\hline 
\# mp & 450 & 230 & {\bf 222} & 235 & 260 & 1\,000 \\
speedup & 2.2 & 4.3 & {\bf 4.5} & 4.3 & 3.8 & 1.0 \\
efficiency & {\bf 22\,\%} & {\bf 22\,\%} & 15\,\% & 11\,\% & 8\,\% & 100\,\% \\
\hline
\end{tabular}
\end{center}
\caption{\label{tab:csk} Errors $|c_s^{(k), {\sfrei \text{fine}}}(T_{\text{end}})-c_s^{*}(T_{\text{end}})|$ for $P=10,\ldots, 50$ for the parareal algorithm (\cref{alg:parareal}) in the first numerical example (ODE growth). {\reftwo For comparison, the reference value $c_s^{*}(T_{\text{end}})=0.63831273...$ resulting from a serial time-stepping with the same fine-scale time-step size $\delta t$ is taken.} The time measure in terms of the number of serial micro problems (\# mp) as well as speedup and efficiency compared to the reference computation (right column) are shown; the stopping criterion $|c_s^{(k+1), {\sfrei \text{fine}}}(T_{\text{end}}) - c_s^{(k), {\sfrei \text{fine}}}(T_{\text{end}})| < \epsilon_{\text{par}}=10^{-3}$ is used; best numbers marked in \textbf{bold face}.}
\end{table}

In~\cref{tab:csk}, we show the deviations in the foam cell concentration $|c_s^{(k), {\sfrei \text{fine}}}(T_{\text{end}})-c_s^*(T_{\text{end}})|$ at final time $T_{\text{end}}$ after each iteration of the parareal algorithm for $P=N_c = 10, 20, 30, 40$ and $50$ processes. 
We observe that the number of iterations $k_{\text{par}}$ decreases from four to three for $P\geq 20$. This is due to the fact that the coarse problem is solved with the smaller time-step size {\sfrei $\delta T = \lceil \frac{T_\text{end}}{P}\rceil$}, which makes the coarse problem more expensive but also more accurate. In particular, we obtain 
a better approximation for the coarse values $\overline{c}_s^{(k)}(T_p)$ that are used as initial values in the fine problems in the next iterate. {\reftwo We can also observe that, for fixed $P$, the error decreases at least by a factor $\frac{1}{P}$ in each parareal iteration. This is in agreement with~\cref{lem.conv_parareal}, which predicts (for a simpler model problem) a reduction factor $c\cdot \delta T = c\cdot \lceil \frac{T_\text{end}}{P}\rceil$ for some constant $c>0$, compared to the previous iterate.}

Since the number of parareal iterations $k_{\text{par}}$ is constant for all cases with $P\geq 20$, we get the lowest computational cost for $P=30$, which is close to $\sqrt{N_f}\approx 31.6$; cf. the discussion 
in~\cref{sec:compcost}. For $P=30$, 
$3\cdot 34=102$ micro problems need to solved 
on each process (step II.a) and $4\cdot 30 = 120$ micro problems within the coarse propagators (steps I and II.b), {\reftwo i.e., $222$ micro problems in total}. {\reftwo Compared to a serial time-stepping scheme with the same fine-scale time-step size, which requires the solution of $N_f=1\,000$ serial micro problems, this corresponds to a speed-up of $\frac{1\,000}{222} \approx 4.5$; {\reftwo see also~\cref{eq:speedup}}.} 
As is usual for two-level methods, for larger  numbers of processes ($P$), the cost of the coarse problems gets dominant. For $P=40$, we solve, for example, $3\cdot 25=75$ micro problems in the fine propagator compared to $4\cdot 40=160$ micro problems for the coarse propagator; this results in a total of $235$ micro problems and a speedup of $4.3$. 

\begin{table}[t]
\begin{center} \small
\begin{tabular}{|l|rrrrr|r|}
\hline
k &$P=10$ &$P=20$ &$P=30$ &$P=40$ &$P=50$ & ref. (serial) \\
\hline                                           
1 &$2.59\cdot 10^{-3}$ &$5.97\cdot 10^{-4}$ &$3.29\cdot 10^{-4}$ &$1.42\cdot 10^{-4}$ &$8.95 \cdot 10^{-5}$ &{\reftwo 0} \\ 
2 &$2.81\cdot 10^{-3}$ &$2.04\cdot 10^{-5}$ &$6.59\cdot 10^{-6}$ &$2.36\cdot 10^{-6}$ &$8.43\cdot 10^{-7}$ &- \\
3 &$6.37\cdot 10^{-6}$ &- &- &- &- &-\\
4 &$1.73\cdot 10^{-7}$ &- &- &- &- &-\\
\hline 
\hline
\# mp &450 &160 &{\bf 158} &170 &190 &1\,000 \\
speedup & 2.2 & 6.25 & {\bf 6.3} & 5.9 & 5.3 & 1.0 \\
efficiency & 22\,\% & {\bf 31\,\%} & 21\,\% & 15\,\% & 11\,\% & 100\,\% \\
\hline
\end{tabular}
\end{center}
\caption{\label{tab:csoverline} Errors $|\overline{c}_s^{(k)}(T_{\text{end}})-c_s^{*}(T_{\text{end}})|$ for $P=10, \ldots, 50$ for the parareal algorithm (\cref{alg:parareal}) in the first numerical example (ODE growth).
The time measure in terms of the number of serial micro problems (\# mp) as well as speedup and efficiency compared to the reference computation (right column) are shown; the stopping criterion $|\overline{c}_s^{(k+1)}(T_{\text{end}}) - \overline{c}_s^{(k)}(T_{\text{end}})| < \epsilon_{\text{par}}=10^{-3}$ is used; best numbers marked in \textbf{bold face}.
}
\end{table}

{\reftwo In order to relate the speed-up to the computational resources used, we define the computational efficiency by the ratio ${\displaystyle \frac{\text{speedup}}{P}}$ as a second measure to quantify the results.
By definition, the serial computation always has the best efficiency,} as it converges within one sweep through all time steps. Even if the parareal algorithm would also converge within one iteration, the parallel computations come with additional costs since they require the solution of the coarse problem. 

In fact,
using the serial computation as a reference in the comparison in~\cref{tab:csk}, the efficiency is best ($22\,\%$) for $P=10$ or $P=20$; it deteriorates for larger numbers of processes due to the increasing cost of the coarse problems. The low efficiency is a typical observation for parallel time integration methods, such as the parareal method and can be explained as follows: For $P=10$, for instance, 
4 parareal iterations are required and each iteration requires the solution of all 1\,000 micro problems on the fine scale on the full time interval $[0,T_{\text{end}}]$. Hence, a total of 4\,000 micro problems has to be solved. Assuming that the time to solve a micro problem is approximately constant,
we can not expect efficiencies above 25\%. Due to the additional cost of the coarse propagator, the efficiency is reduced to 22\%. However, since we can always compute 10 micro problems on the fine scale at the same time, we can still obtain a speedup of 2.2. 

To increase the efficiency, we could take the value $\overline{c}_s^{(k)}(T_{\text{end}})$ computed in~\cref{eq:parareal} instead of $c_s^{(k), {\sfrei \text{fine}}}(T_{\text{end}})$ in the stopping criterion
\begin{align}\label{stoppingcrit}
|\overline{c}_s^{(k+1)}(T_{\text{end}}) - \overline{c}_s^{(k)}(T_{\text{end}})| < \epsilon_{\text{par}}.
\end{align}
In particular, at 
the end time $T_{\text{end}}$ we expect that this has a higher accuracy compared to the fine-scale variable $c_s^{(k), {\sfrei \text{fine}}}(T_{\text{end}})$.
However, the values $\overline{c}_s^{(k+1)}(t_i)$ are only available at the coarse grid points $t_i\in \{T_p\}_{p=1}^P$. If one is interested in 
foam cell concentrations at intermediate time steps, a stopping criterion based on $c_s^{(k), {\sfrei \text{fine}}}(t_i)$
needs to be used, as in~\cref{eq:terminal}.

The errors concerning $\overline{c}_s^{(k)}(T)$ and the number of micro problems to be solved using the stopping criterion~\cref{stoppingcrit} are given in~\cref{tab:csoverline}. We observe that for $P\geq 20$ this stopping criterion was satisfied already after $k_{\text{par}}=2$ iterations. For $P=30$ we get again the lowest computational cost with a speed-up of around $6.33$ compared to a serial computation. Furthermore, the best efficiency is obtained here with $P=20$, which is due to the lower iteration count compared to the case $P=10$.

\begin{figure}[t]
\centering
\begin{minipage}{0.52\textwidth}
\hspace*{-0.5cm}\includegraphics[width=\textwidth]{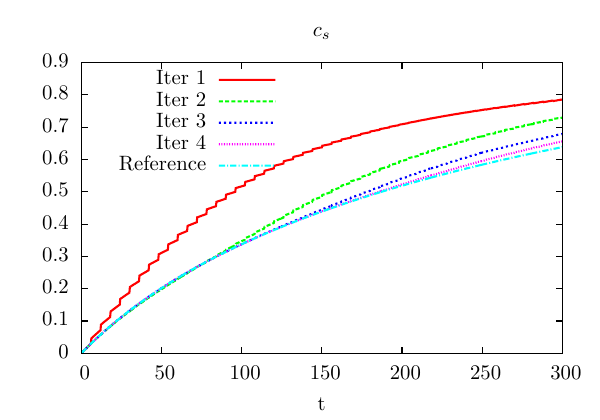}
\end{minipage}\hspace*{-0.3cm}
\begin{minipage}{0.52\textwidth}
\hspace*{-0.5cm}\includegraphics[width=\textwidth]{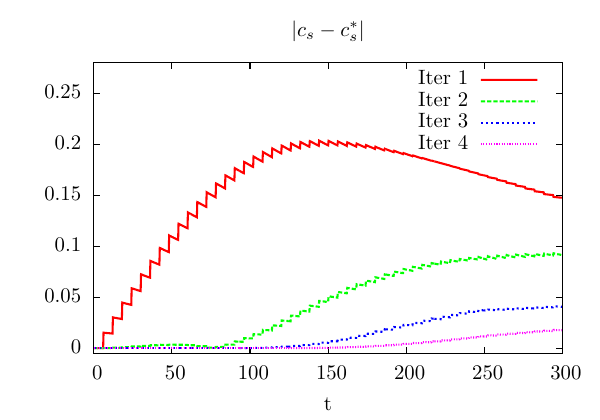}
\end{minipage}\\
\begin{minipage}{0.45\textwidth} \small
\begin{tabular}{|l|rr|r|} 
\hline
k &$P=30$ &$P=50$ 
& ref. (serial) \\
\hline
1 &$1.66\cdot 10^{-1}$ &$1.42\cdot 10^{-1}$
  &{\reftwo 0}\\
2 &$9.34\cdot 10^{-2}$ &$9.34\cdot 10^{-2}$
 &-\\
3 &$4.86\cdot 10^{-2}$ &$4.06\cdot 10^{-2}$
&- \\
4 &$1.84\cdot 10^{-2}$ &$1.76\cdot 10^{-2}$
&- \\
5 &$5.41\cdot 10^{-3}$ &$5.71\cdot 10^{-3}$
&-\\
6 &$1.28\cdot 10^{-3}$ &$1.47\cdot 10^{-3}$
&- \\
7 &$2.47\cdot 10^{-4}$ &$3.14\cdot 10^{-4}$ 
&- \\
8 &$4.79\cdot 10^{-5}$ &$6.48\cdot 10^{-5}$ 
&-\\
\hline \hline 
\# mp &272 & {\bf 160} 
&1\,000 \\
speedup & 3.7 & {\bf 6.3} 
&1.0 \\
eff. & 12\,\% & {\bf 13\,\%} 
& 100\,\% \\
\hline
\end{tabular}
\end{minipage}
\caption{\label{fig:pararealmacro} Behavior of the iterates of the parareal algorithm, when the stationary FSI problem (\cref{long-scale}) is used in Step (I) and (II.b)(i) of~\cref{alg:parareal}. A description of the computation is given in~\cref{sec:pararealmacro}. Convergence to the reference solution is much slower compared to standard parareal (see~\cref{fig:parallel}). \textbf{Top}: {\reftwo Evolution of $c_s^{(k), {\text{fine}}}(t)$ and error $|c_s^{(k), \text{fine}}(t) - c_s^*(t)|$ of the} foam cell concentration over time $t$ for $P=50$ and $k=1,2,3,4$. \textbf{Right}: Errors $|c_s^{(k), {\sfrei \text{fine}}}(T_{\text{end}})-c_s^{*}(T_{\text{end}})|$ after each iteration for $P=30$ and $P=50$; total number of micro problems (\# mp) solved in parallel as well as speedup and efficiency compared to the reference computation (right column); best numbers marked in \textbf{bold face}.
}
\end{figure}

\subsection{Variants with cheaper coarse-scale computations}
\label{sec:pararealmacro}

In the parareal algorithm introduced above, a further 
improvement of the computational cost is not possible for $P \gtrsim \sqrt{N_f}$ processes {\reftwo in the case $P=N_c$} due to the (increasing) cost of the coarse-scale propagators. In fact, these get dominant compared to the fine-scale contributions for $P\gtrsim \sqrt{N_f}$. In this section, we will discuss approximate coarse-scale propagators, where no additional computations of micro problems are needed; hence, they are computationally cheaper.

\subsubsection{Heuristic averaging of the FSI problems}

As a first variant, 
we use the heuristic averaging approach mentioned in the beginning of~\cref{sec:two_scale} for the coarse-scale propagation in steps (I) and (II.b)(i) of~\cref{alg:parareal}.
The solution of the stationary FSI problem in~\cref{long-scale} is much cheaper compared to $\geq 100$ time steps of a non-stationary FSI problem (approx.\,by a factor 100), and thus, its computational cost is 
neglected in the following discussion. 
Hence, the computational cost is reduced from $k_{\text{par}}\cdot \lceil N_f/P \rceil + (k_{\text{par}}+1) \cdot P$ to $k_\text{par} \cdot \lceil N_f/P \rceil$ micro problems. In~\cref{fig:pararealmacro}, we {\reftwo illustrate the error in} the foam cell concentration for $P=50$ (left), and give the {\reftwo error $|c_s^{(k), {\sfrei \text{fine}}}(T_{\text{end}})-c_s^*(T_{\text{end}})|$ in each iteration $k$ for $P=30$ and $P=50$ (right)}. We observe a much slower convergence compared to the standard parareal algorithm used above. This could already be expected from the results in~\cref{fig:NumEx}, where we saw that the heuristic averaging is not a good approximation. Using the stationary problem for the coarse propagator, $8$ parareal iterations are necessary until the stopping criterion is satisfied for $P=30$ and $P=50$ compared to {\sfrei $3$ iterations} in the standard parareal algorithm.

In terms of the 
computating time, 
$8\cdot 34=272$ and $8\cdot 20=160$ micro problems are necessary for the cases $P=30$ and $P=50$, respectively. While for $P=30$ this is worse compared to the standard parareal algorithm, this is an improvement of $67\%$ for $P=50$; cf.~\cref{tab:csk}.
As for the standard parareal algorithm, convergence is obtained faster
if we consider the stopping criterion in~\cref{stoppingcrit} based on the parareal iterates $\overline{c}_s^{(k)}(T_{\text{end}})$; however, $7$ iterations are still needed.
The resulting computational cost is $7\cdot 34=238$ micro problems for $P=30$ and $7\cdot 20=140$ micro problems for $P=50$. Compared to 158 micro problems for $P=30$ and 190 micro problems for $P=50$ in the standard parareal algorithm (cf.~\cref{tab:csoverline}), we see only a small but no significant improvement. We will thus consider another, more promising approach to compute the required growth values $\gamma(\sigmat_f^{WS})$ for the coarse propagator in the following subsection.

\begin{remark}
 As mentioned before, parallelization in space is generally more efficient than parallelization in time. Hence, complex three-dimensional FSI simulations already have a significant potential for parallelization. Additionally using temporal parallelization with large values of $P$, such as $P>50$, may only be reasonable on very large supercomputers.
\end{remark}

\subsubsection{Re-usage of computed growth values}
\label{sec:reusage}

As a second approach we will consider 
the re-use of growth values $\overline{\gamma}(\sigmat_f^{WS}(t_{p,i}))$ computed in the fine-scale propagator on the coarse scale. 
Therefore, we store all values $\overline{\gamma}_{p\cdot n_p + i}:=\overline{\gamma}(\sigmat_f^{WS}(t_{p,i}))$ computed on the fine scale on all processes $p=1,...,P$ 
for all time-steps $i=1,...,n_p$; see step 1.b of~\cref{twoscale}. {\reftwo These can be used in the coarse propagator (step II.b) in the same parareal iteration instead of computing new growth values there.} 

{\reftwo We introduce the following notations for the coarse resp.~fine propagators that start from $c_{s,n-1}$ using certain growth values $\gamma(T_{n-1})$ (that might now differ from $\overline{\gamma}(\sigma_f^{\text{WS}}(c_{s,n-1}))$): 
\begin{align*}
c_{s,n}^{\text{coarse}} = {\cal C}(I_n,c_{s,n-1}, \gamma(T_{n-1})),\quad
c_{s,n}^{\text{fine}} = {\cal F}(I_n,c_{s,n-1}, \gamma(T_{n-1})),
\end{align*}
where $I_n = [t_{n-1}, t_n]$.
After an initialization step, the modified parareal iteration is defined by the following formula for $n=1,...,P$ and $k\geq 1$:
\begin{align}\label{eq:parareal_app}
\overline{c}_{s,n}^{(k)} = {\cal C} (I_n,\overline{c}_{s,{n-1}}^{(k)}, \bar{\gamma}(c_{s,n-1}^{(k), \text{fine}})) +{\cal F}(I_n, \overline{c}_{s,n-1}^{(k-1)}, \bar{\gamma}(\overline{c}_{s,n-1}^{(k-1)})) - {\cal C}(I_n,\overline{c}_{s,n-1}^{(k-1)}, \bar{\gamma}(c_{s,n-1}^{(k-1),\text{fine}})),
\end{align}
where
\begin{align}\label{Defcsfine}
c_{s,n-1}^{(k),\text{fine}} := {\cal F}(I_{n-1}, \overline{c}_{s,n-2}^{(k-1)}, \bar{\gamma}(\overline{c}_{s,n-2}^{(k-1)})),
\end{align} 
and $c_{s,n}^{(0), \text{fine}}=\overline{c}_{s,n}^{(0)}$. Moreover, we set $ c_{s,0}^{(k), \text{fine}} = \overline{c}_{s,0}^{(k)} =0$ for all $k$.}


As no new micro problems need to be solved and approximations of the growth values $\overline{\gamma}_j$ are available for all fine time steps $j=1,...,N_f$, the coarse propagator can now {\sfrei even} use the fine-scale time-step $\delta t$. The only additional cost is to advance the foam cell concentration by~\cref{reaction} resp.~\cref{eq:parareal} {\sfrei on the coarse level}. This cost is clearly negligible for the ODE growth model in~\cref{reaction}. The additional cost in case of the PDE model~\cref{eq:parareal} will be discussed in~\cref{sec:reaction-diffusion}. The resulting algorithm is given as~\cref{alg:para}.

\begin{algorithm}[t]
	\caption{Parallel Time-Stepping with Re-Usage of Growth Values\label{alg:para}}
	\begin{enumerate}
		\item[(I)] \textbf{Initialization}: Compute $\big\{(\overline{c}_s^{(0)}(T_p), \overline{\wt}^{(0)}(T_p)\big\}_{p=1}^P$ by means of~\cref{twoscale} 
		with a coarse macro time-step size $\delta T=T_{p+1}- T_p$ on the master process. Set $k\gets 0$\\
		\item[(II)]\While{$|c_s^{(k+1) ,{\sfrei \text{fine}}}({\sfrei T_{\text{end}}}) - c_s^{(k), {\sfrei \text{fine}}}({\sfrei T_{\text{end}}})| > \epsilon_{\text{par}}$}{
			\begin{itemize}
				\item [(II.a)] \textbf{Fine problem}: \\
					\For{$p=1,...,P$}{ 
				\begin{itemize}
					\item[(i)] Initialize $c_s^{(k+1), {\sfrei \text{fine}}}(T_p)=\overline{c}_s^{(k)}(T_p), \,\wt^{(k+1)}(T_p)  = \wt^{(k)}(t_{p-1, N_{p}}))$
					\item[(ii)] Compute $\big\{(c_s^{(k+1), {\sfrei \text{fine}}}(t_{p,q}), \wt^{(k+1)}(t_{p,q})\big\}_{q=1}^{n_p}$ by~\cref{twoscale} with fine time step $\delta t$
					\item[(iii)] Store the resulting growth functions {\reftwo $\overline{\gamma}_{p\cdot n_p+q}:=\overline{\gamma}(\sigmat_f^{WS})(t_{p,q})$ at all fine points \\$t_{p,q}, \,q=1,...,n_p$} as well as $\wt(t_{p,n_p})$ at the last time-step.
				\end{itemize}  
			}
				\item[(II.b)] \textbf{Coarse problem}:
				\begin{itemize}
					\item[(i)] \For{$j=1,...,N_f$}{
						\hspace*{2cm}Compute $\overline{c}_s^{(k+1)}(t_j)$ by 
						advancing the ODE~\cref{reaction} resp.\, solving the PDE\\ \hspace*{1.95cm} in~\cref{PDEmodel}  using the growth values $\overline{\gamma}_j$ computed in (II.a).
					}
				\end{itemize}
		\end{itemize}}
	\end{enumerate}
\end{algorithm}

The growth values employed for the re-usage are available for all fine-scale time steps, and hence, we are able to perform the coarse propagation efficiently on the fine scale; cf.~the discussion on the computational cost later in this section as well as for the PDE growth model 
in~\cref{sec:runtimes}. Nonetheless, the expected accuracy of the re-usage approach is still lower compared to the standard parareal algorithm. This is because the growth values on the fine scale have been computed using {\sfrei $c_{s,n-1}^{(k),\text{fine}}$, which depends on the previous iterate $\overline{c}_{s,n-2}^{(k-1)}$ (see~\cref{Defcsfine})}, whereas the coarse propagator in the standard parareal algorithm already uses the more accurate {\sfrei $\overline{c}_{s,n-1}^{(k)}$} from the current iteration.

\paragraph{Computational costs}

Only the coarse step in the initialization (I) has to be carried out without re-usage
and comes with a computational cost of {\reftwo $N_c$} micro problems. On the other hand,
step (II.b) does not require the solution of any micro problems. In terms of micro problems to be solved, the computational cost of step (II.a) is the same as in
the standard parareal algorithm (\cref{alg:parareal}).
Altogether, the number of micro problems to be solved in $k_{\text{par}}$ iterations of~\cref{alg:para} is
\begin{align*}
k_{\text{par}}\cdot \lceil N_f/P\rceil + {\reftwo N_c}.
\end{align*}
{\reftwo Again, we obtain a saturation of the cost for large $P$, but the cost for large $P$	
is by a factor $k_{\text{par}}$ smaller. Moreover, for $N_c=P$,} the choice $P\approx \sqrt{k_{\text{par}} N_f}$ would be optimal if the number of iterations $k_{\text{par}}$ was independent of the number of processes $P$.

{\reftwo The speed-up compared to a serial computation is given by
\begin{align}\label{eq:speedup_reusage}
\text{speedup}(P) = \frac{N_f}{k_{\text{par}}\cdot \lceil N_f/P \rceil  + N_c} \approx \frac{1}{\frac{k_{\text{par}}}{P}  + \frac{N_c}{N_f}},
\end{align}
assuming, again, a perfect load balancing of the micro problems.}

{\reftwo The communication cost of the approach strongly depends on the employed parallelization scheme; cf.~\cref{sec:parallelization}. In the distributed parallelization approach, the growth values $\overline{\gamma}(\sigmat_f^{WS}(t_{p,i}))$, which are needed in the coarse propagator, do not have to be communicated. This is because the same intervals in fine and coarse problems are computed on the processes. In the master-slave approach, however, the values have to be communicated from each slave process to the master.}
Using the ODE growth model, these are $N_f$ scalar values in total (in our example $N_f=1\,000$). {\reftwo Even for the master-slave communication scheme, it} is thus reasonable to assume that the computational cost of communication is still negligible compared to the solutions of the micro problems. In the case of a PDE growth model, the additional cost for communication will be discussed in~\cref{sec:reaction-diffusion}.

{\reftwo \paragraph{Theoretical convergence analysis}

We extend the convergence results discussed in~\cref{subsec.conv} for the model problem in~\cref{eq:EulerODEcont} for the re-usage algorithm in~\cref{eq:parareal_app} under the assumptions made in~\cref{subsec.conv}.

The following result is a direct consequence of~\cref{lem.conv_reusage}, which is shown in the appendix.

\begin{theorem}\label{theo.conv_reusage}
Let $c_s^{*}\in C^1(0,T_\text{end})$ be the exact solution of~\cref{eq:EulerODEcont} and let $\{\overline{c}_s^{(k)}(T_n)\}_{n=1}^{P}$ be the $k$-th iterate of the re-usage algorithm in~\cref{eq:parareal_app} using the forward Euler method in~\cref{eq:EulerODE2}. Under the assumptions made in~\cref{subsec.conv}, it holds for the error $ e_n^{(k)} = |\overline{c}_s^{(k)}(T_n)-c_s^{*}(T_n)|$ and $k\in\mathbb{N}_0, n\in \{1,\ldots,P\}$ that
\begin{align}
e_n^{(k)} 
&\leq L \delta T  \max\left\{1,T_{n-\frac{k}{2}}\right\}^k \frac{c^k \beta^{n-k}}{\lceil \nicefrac{k}{2} \rceil !} \max_{t\in [0,T_\text{end}]} |\partial_t c_s^{*}(t)|\label{eq:est_reusage}.
\end{align}
with a constant $c>0$ (see~\cref{lem.conv_reusage})
and $\beta=1+L\delta T$. This implies $e_n^{(k)}\to 0$ for $k\to\infty$.
\end{theorem}
The convergence $e_n^{(k)}\to 0$ for $k\to\infty$ follows due to the faculty ${\lceil \nicefrac{k}{2} \rceil !}$ in the denominator. Compared to~\cref{lem.conv_parareal}, the estimated convergence in $k$ is much slower. However, we will observe in the numerical examples below that a good accuracy might still be reached within few iterations.
}

\begin{figure}[t]
\center
\hspace*{-0.4cm}\includegraphics[width=0.52\textwidth]{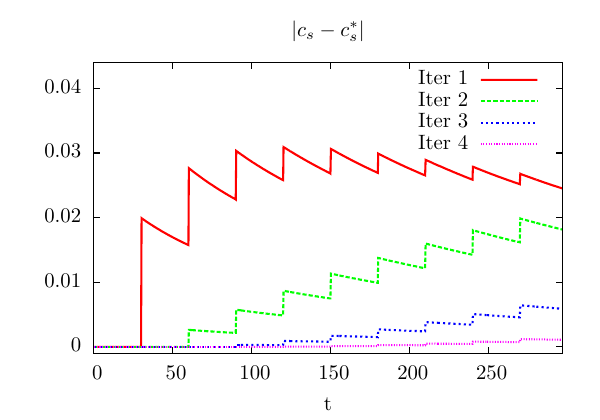}\hspace*{-0.5cm}
\includegraphics[width=0.52\textwidth]{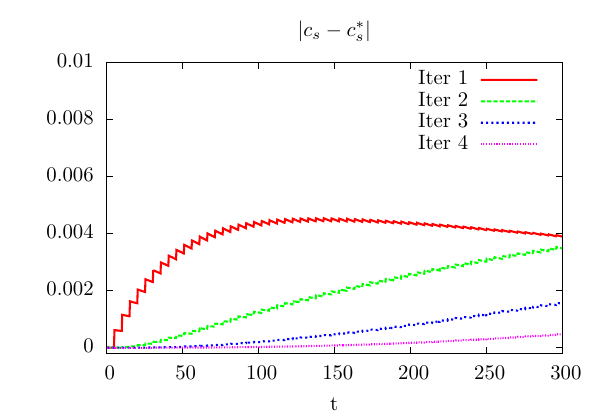}
\caption{\label{fig:reusage} {\reftwo Error in the concentration $|c_s^{(k),\text{fine}}(t)-c_s^*(t)|$} over time $t$ in the first 3 iterates of the parareal algorithm with re-usage of growth values (\cref{alg:para}) for $P=10$ (left) and $P=60$ (right) and the first numerical example (ODE growth model).
}
\end{figure}

\begin{table}[t]
\begin{center}
\centering
\setlength{\tabcolsep}{4pt}\renewcommand{\arraystretch}{1.1} \small
\begin{tabular}{|l|rrrrrrr|r|}
\hline
k &$P=10$ &$P=20$ &$P=30$ &$P=40$ &$P=50$ &$P=60$ &$P=70$ & ref. \\
\hline                                            
1 &\small $2.56\cdot 10^{-2}$ &\small $1.10\cdot 10^{-2}$ &\small $7.48\cdot 10^{-3}$ &\small $5.46\cdot 10^{-3}$ &\small $4.23\cdot 10^{-3}$ &\small$3.63\cdot 10^{-3}$ &\small $3.12\cdot 10^{-3}$ &\small {\reftwo 0}\\
2 &\small $7.78\cdot 10^{-3}$ &\small $4.48\cdot 10^{-3}$ &\small $3.16\cdot 10^{-3}$ &\small $2.36\cdot 10^{-3}$ &\small $1.99\cdot 10^{-3}$ &\small $1.61\cdot 10^{-3}$ &\small $1.41\cdot 10^{-3}$ &\small-\\
3&\small $1.73\cdot 10^{-3}$ &\small $1.20\cdot 10^{-3}$ &\small $9.02\cdot 10^{-4}$ 
&\small $6.92\cdot 10^{-4}$ &\small $6.34\cdot 10^{-4}$ &\small $5.14\cdot 10^{-4}$ &\small $4.32\cdot 10^{-4}$ &\small-\\
4 &\small $2.32\cdot 10^{-4}$ &\small $2.28\cdot 10^{-4}$ &\small $1.87\cdot 10^{-4}$ &\small $1.31\cdot 10^{-4}$ &\small $1.26\cdot 10^{-4}$ &\small $1.10\cdot 10^{-4}$ &\small $9.98\cdot 10^{-5}$  &\small-\\
5 &\small $2.19\cdot 10^{-5}$ &\small $4.71\cdot 10^{-5}$ &\small $3.05\cdot 10^{-5}$ &\small - &\small - &\small - &\small - &\small -\\
\hline
\hline 
\# mp 	&510 &270 &200 &140 &130 &{\bf 128} &130 &1\,000 \\
speedup &2.0 &3.7 &5.0 &7.1 &7.7 & {\bf 7.8} &7.7 &1.0 \\
eff. 	&{\bf 20\,\%} &19\,\% &17\,\% &18\,\% &15\,\% &13\,\% &11\,\% &100\,\%  \\
\hline
\end{tabular}
\end{center}
\caption{\label{tab:para2} 
Errors $|c_s^{(k), {\sfrei \text{fine}}}(T_{\text{end}})-c_s^{*}(T_{\text{end}})|$ for $P=10, \ldots, 70$ for the parareal~\cref{alg:para} (re-use of growth values) and the first numerical example (ODE growth model). Value {\sfrei$c_s^*(T_{\text{end}})$} of a {\sfrei serial} reference computation for comparison. 
The time measure in terms of the number of micro problems (\# mp) as well as speedup and efficiency compared to the reference computation (right column) are shown; the stopping criterion $|c_s^{(k+1), {\sfrei \text{fine}}}(T_{\text{end}}) - c_s^{(k), {\sfrei \text{fine}}}(T_{\text{end}})| < \epsilon_{\text{par}}=10^{-3}$ is used; best numbers marked in \textbf{bold face}.
}
\end{table}

\paragraph{Numerical results}

In~\cref{fig:reusage}, the evolution of the {\reftwo error in the} concentration variable $c_s^{(k)}$ is plotted over time for {\reftwo $P=N_c=10$} and {\reftwo $P=N_c=60$}. We observe convergence to the reference values $c_s^*$ in both cases. Compared to the results in~\cref{fig:parallel} for the classical parareal algorithm (\cref{alg:parareal}), {\reftwo the convergence is significantly slower}. Nevertheless, the stopping criterion is already fulfilled after 4-5 iterations (see~\cref{tab:para2}), which is only 1-2 iterations more than in the standard parareal algorithm (see~\cref{tab:csk}).

{\reftwo While we had observed reduction factors in the order of $\delta T = \lceil \frac{T_{\text{end}}}{P}\rceil$ between two consecutive parareal iterations for the standard parareal algorithm in~\cref{tab:csk}, the reduction factors are slightly worse in~\cref{tab:para2}. They are, however, in all cases (besides the last value for $P=20$) bounded above by $\frac{1}{k}$, where $k$ is the parareal iterate. This indicates a convergence of order ${\cal O}(\frac{1}{k!})$, which is related to the factor ${\cal O}(\frac{1}{\lceil \frac{k}{2}\rceil !})$ in~\cref{eq:est_reusage} in~\cref{theo.conv_reusage}.
}
{\reftwo Moreover, the absolute values of the error in each parareal iteration are significantly smaller for $P=60$ compared to the case of $P=10$ (Note the different scaling in the vertical axis)}.

In~\cref{tab:para2}, we compare the number of micro problems needed for {\reftwo $P=N_c=10$} to $70$. 
The optimal number of processes is $P=\sqrt{k_{\text{par}} N_f} \approx 63.25$, which is twice as many processes compared to the parareal algorithm in the previous section. The minimal cost in~\cref{tab:para2} is 128 micro problems for $P=60$, compared to 222 for standard parareal. Compared to a serial time-stepping scheme, we obtain a maximum estimated  
speed-up of 7.8. Moreover, the efficiency is also significantly improved for larger numbers of processes compared to the standard approach; cf.~\cref{tab:csk,tab:para2}.


\section{Plaque growth problem with a distributed foam cell concentration} \label{sec:reaction-diffusion}


We consider the PDE growth model introduced in~\cref{sec:PDE}, where the concentration of foam cells $\hat{c}_s = \hat{c}_s(\hat{x},t)$ is now governed by the non-stationary reaction-diffusion problem in~\cref{PDEmodel}.
For time discretization, we use a linearized variant of the backward Euler scheme. 
Starting from $c_{s,0} = 0$, the standard backward Euler scheme writes for $l=1,...,N_f$: \\ \textit{Find ${\sfrei c_{s,l+1}} \in {\cal C}:= H^1_0({\hat{\SO}},\hat{\Gamma}_s)$ such that}
{\sfrei 
\begin{align}
\begin{split}\label{eq:bwEuler}
\frac{1}{\delta t} \left(\hat{c}_{s,l+1} - \hat{c}_{s,l},\hat\varphi\right)_{\hat{\SO}}   - \left(D_s \hat{\nabla} \hat{c}_{s,l+1}, \hat{\nabla} \varphi\right)_{\hat{\SO}} 
+ R_s \left(\hat{c}_{s,l+1}(1-\hat{c}_{s,l+1}), \hat\varphi\right)_{\hat{\SO}}\\
 + \left(\overline{\gamma}(\sigmat_f^{WS}), \hat\varphi\right)_{\hat{\Gamma}} &=0 \quad \forall \hat{\varphi} \in {\cal C}.
 \end{split}
\end{align}}
Since~\cref{eq:bwEuler} is a nonlinear system of equations, several iterations of a nonlinear solver (e.g., Newton's method) would be necessary to solve it. Thus, we consider the following linearized variant, which can be seen as an implicit-explicit (IMEX) scheme (see, e.g.,~\cite{AscherRuuthWetton1995})
{\sfrei 
\begin{align}
\begin{split}\label{imex}
\frac{1}{\delta t} \left(\hat{c}_{s,l+1} - \hat{c}_{s,l},\hat\varphi\right)_{\hat{\SO}}   {\sfrei +} \left(D_s \hat{\nabla} \hat{c}_{s,l+1}, \hat{\nabla} \varphi\right)_{\hat{\SO}} 
+ R_s \theta \left(\hat{c}_{s,l+1}(1-\hat{c}_{s,l}), \hat\varphi\right)_{\hat{\SO}}\\ + R_s (1-\theta) \left(\hat{c}_{s,l}(1-\hat{c}_{s,l+1}), \hat\varphi\right)_{\hat{\SO}} + \left(\overline{\gamma}(\sigmat_f^{WS}), \hat\varphi\right)_{\hat{\Gamma}} &=0,
\end{split}
\end{align}}
where $\theta\in[0,1]$ is a weighting parameter.
In our numerical experiments, we choose
$\theta=0.7$. 
For spatial discretization, we use again $Q_2$ finite elements on the mesh described in~\cref{sec:numex1}. {\reftwo If we assume that evolution of the wall shear stress $\sigmat_f^{WS}$ was given exactly, the following bound could be shown for the discretization error
\begin{align}\label{DiscErr}
\|\hat{c}_s(T_{\text{end}}) - \hat{c}_s^{N_f}\|_{\hat{S}} + \left(\sum_{l=1}^{N_f} \delta T \|\hat{\nabla} (\hat{c}_s(T_l) - \hat{c}_{s,l})\|_{\hat{S}}^2\right)^{1/2} \leq \,c_1 h^2 + c_2 \delta T.
\end{align}
Due to the nonlinear interaction between FSI and reaction-diffusion equation, an analysis of the discretization error of the coupled problem is more involved and not within the scope of this paper. The chosen time-step $\delta T$ is $\frac{1}{1000}$ times the macro time interval length $T_{\text{end}}$, while the horizontal cell size $h$ is a factor $\frac{1}{20}$ of the length of the channel in horizontal direction. Thus, the errors in~\cref{DiscErr} should be roughly equilibrated.
}

\begin{figure}
\begin{center}
\begin{minipage}{0.15\textwidth}
\vspace*{-25mm}
$t = 50$ days
\end{minipage}
\begin{minipage}[t]{0.8\textwidth}
\includegraphics[width=\textwidth]{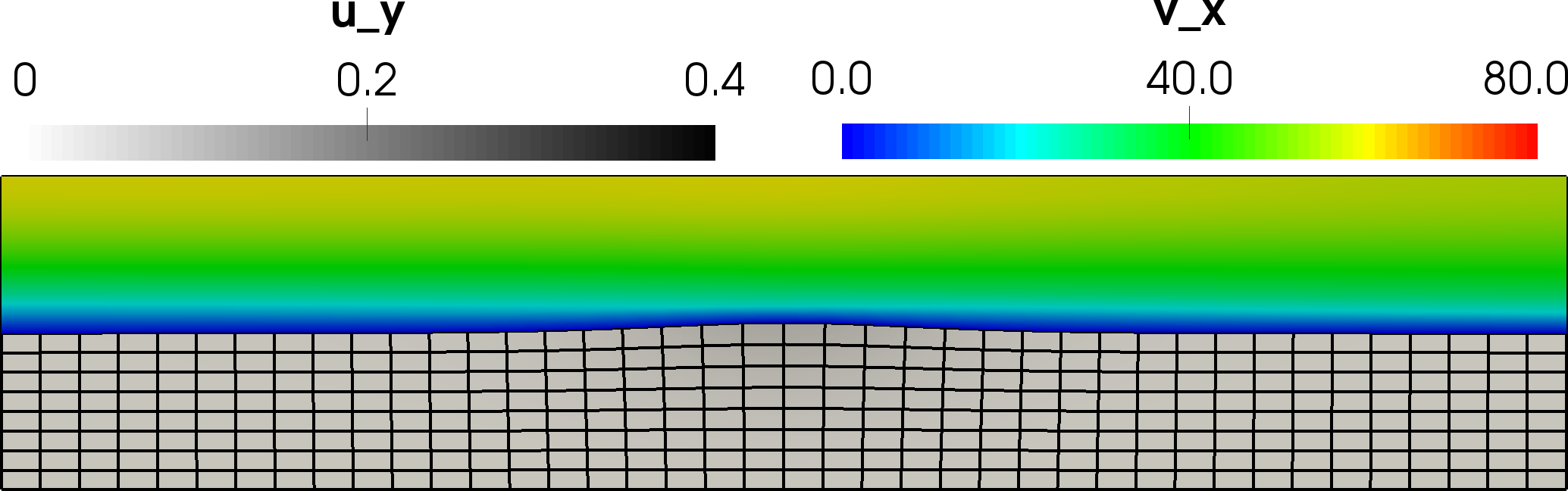}
\end{minipage} \\[1mm]
\begin{minipage}{0.15\textwidth}
\vspace*{-25mm}
$t = 100$ days
\end{minipage}
\begin{minipage}[t]{0.8\textwidth}
\includegraphics[width=\textwidth]{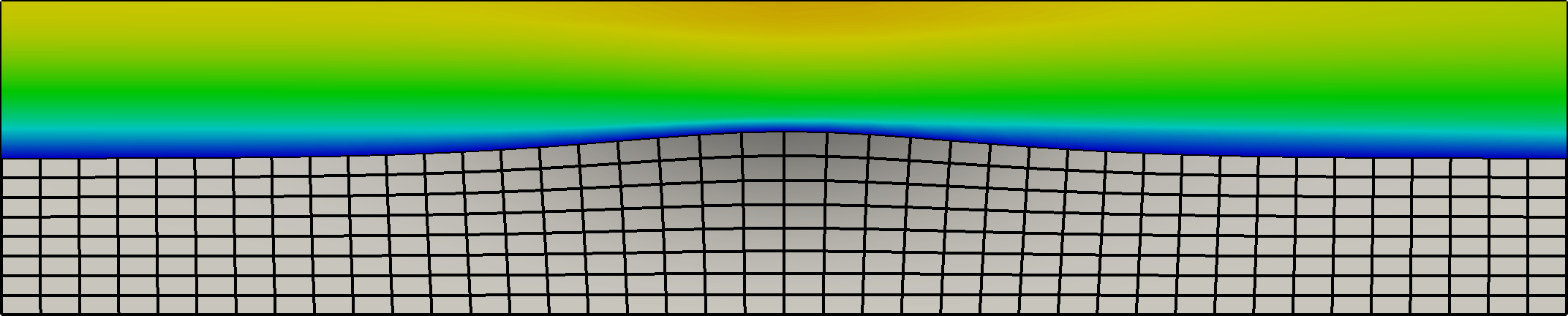} 
\end{minipage} \\[1mm]
\begin{minipage}{0.15\textwidth}
\vspace*{-25mm}
$t = 200$ days
\end{minipage}
\begin{minipage}[t]{0.8\textwidth}
\includegraphics[width=\textwidth]{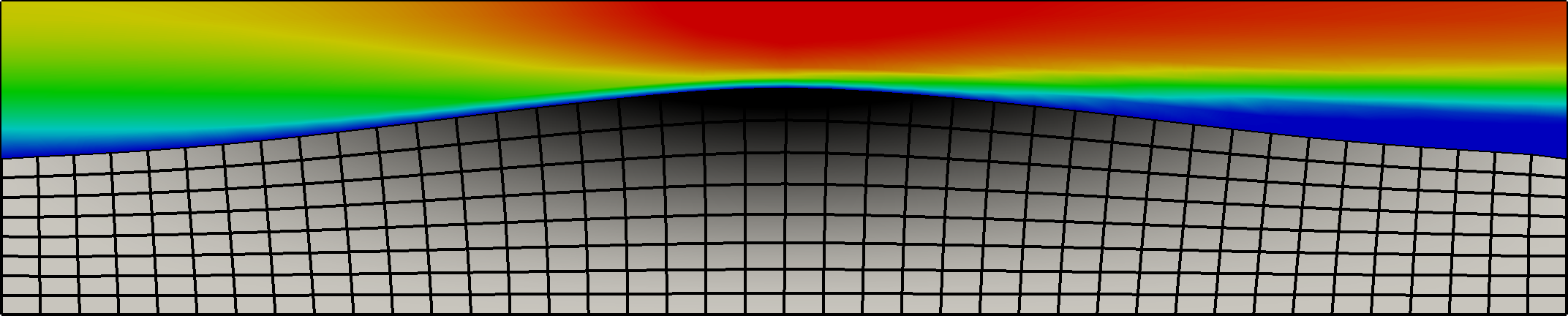}
\end{minipage}
\end{center}
\caption{\label{fig:visu} Visualization of the plaque growth at times $t=50$ days, $t=100$ days and $t=200$ days for the second numerical example (PDE growth model). The horizontal velocity (in cm/s) and the vertical displacement (in cm) are shown on the deformed domain at micro time $\tau=0.5$\,s, i.e., the time of maximum inflow velocity. As the plaque growth evolves, significantly higher velocities arise in the central part.}
\end{figure}

\begin{figure}
\centering
\includegraphics[width=0.7\textwidth]{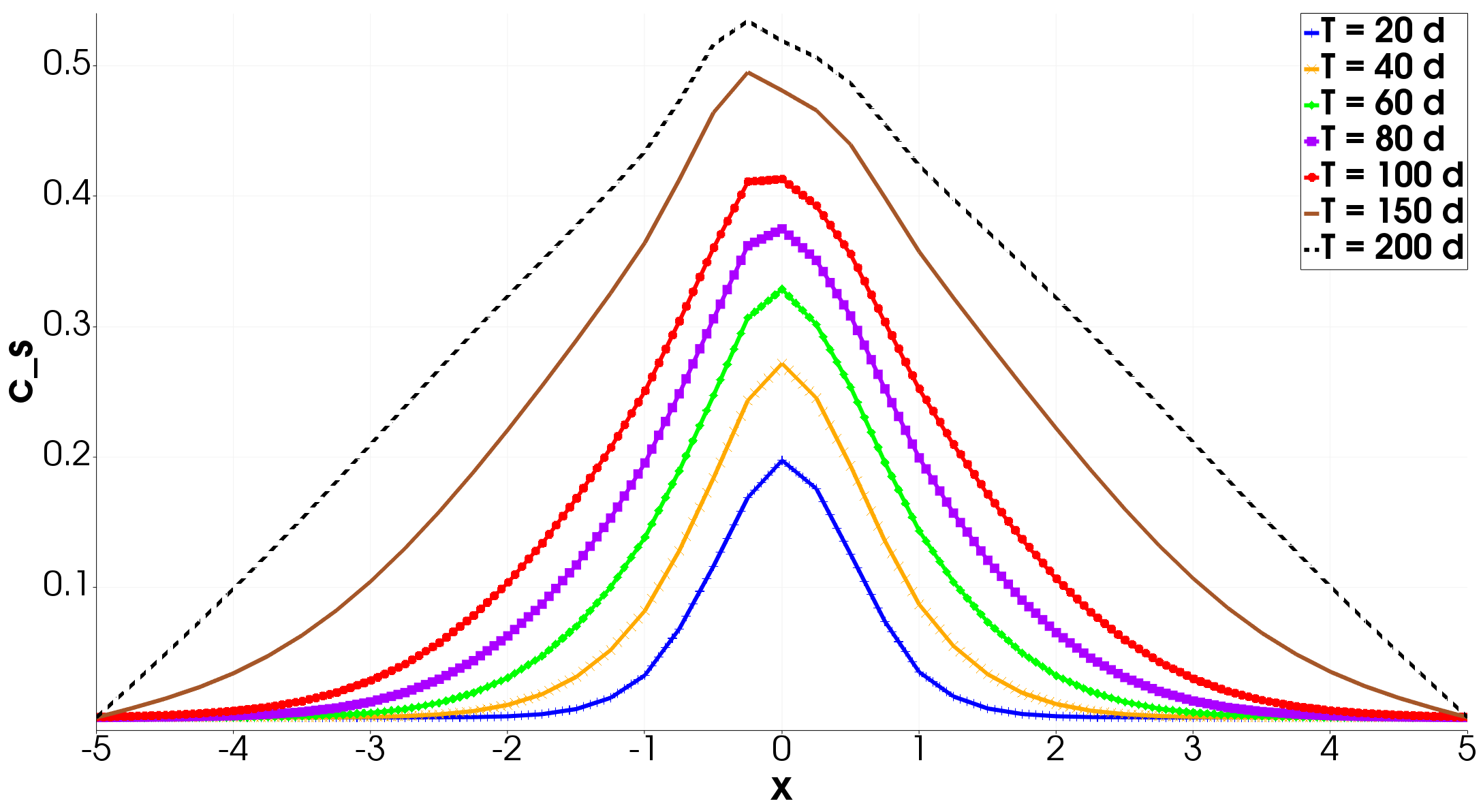}
\caption{\label{fig:results} Foam cell concentration $\hat{c}_s$ at the FSI interface $\hat{\Gamma}$ at different times $t$ in the second numerical example (PDE growth model). In contrast to the ODE model, the concentration is not symmetric around the center anymore for larger times $t$; this is due to the reaction-diffusion equation.}
\end{figure}

The material parameters in the fluid and solid problems are chosen as in the first example. The parameters of the PDE growth model are set to $D_s = 1.2\cdot 10^{-7} \frac{\text{cm}^2}{s}, R_s=5\cdot 10^{-7}\frac{1}{\text{s}}, \alpha=5\cdot 10^{-8}\frac{\text{cm}}{\text{s}}$, and $\sigma_0=30 \frac{\text{g}}{\text{cm s}^2}$. The inflow profile of velocity is chosen as 
\begin{equation*}
  \vt_f^\text{in}(t,x) = 30\begin{pmatrix}
   (1+\sin(\pi \frac{t}{\refone \cal P})^2) (1-x_2^2)\\ 0
  \end{pmatrix}  \text{cm}/\text{s},
\end{equation*}
{\refone where ${\cal P}=1$s}. The end time is $T_{\text{end}}=200$ days, and for the fine time steps, we use 
$\delta t=0.2$ days, such that again $N_f=1\,000$. In~\cref{fig:visu}, we visualize the results for horizontal velocity $v_x$ and {\refone vertical} displacement $u_y$ in the deformed fluid and solid domains in three different time instants, respectively. In~\cref{fig:results} the foam cell concentration $c_s$ on the FSI interface is shown 
for different times $t$. First, we observe a dominant growth in the center of the domain due to the penetration of monocytes and the reaction term in~\cref{PDEmodel}. After $t>100$ the diffusion gets more significant, such that foam cells are distributed over the whole interface.


\begin{figure}[t]
\centering
\hspace*{-0.4cm}\includegraphics[width=0.52\textwidth]{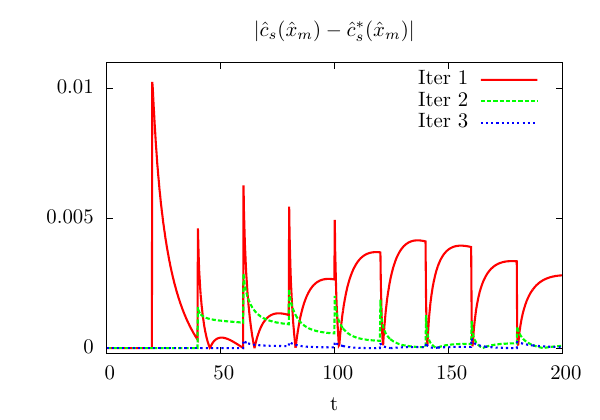}\hspace*{-0.5cm}
\includegraphics[width=0.52\textwidth]{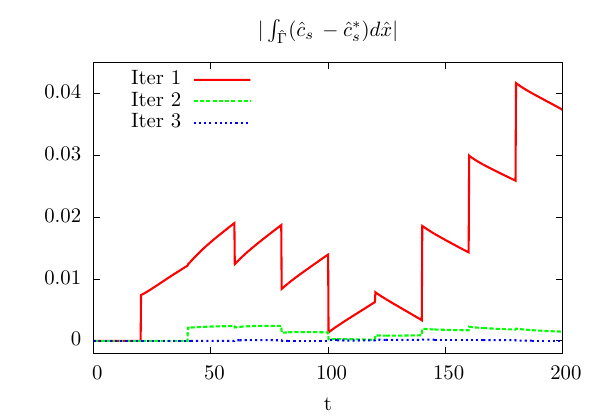}
\caption{{\reftwo Errors} for the second numerical example (PDE growth model) in different iterations of the parareal algorithm for $P=10$; the setting is described in~\cref{sec:reaction-diffusion}. \textbf{Left}: {\reftwo Error in the concentration $|\hat{c}_s^{(k), \text{fine}}(\hat{x}_m,t) - \hat{c}_s^*(\hat{x}_m,t)|$} at the midpoint of $\hat{\Gamma}$ over time; see also~\cref{fig:results} for the temporal evolution of $\hat{c}_s$ at the interface. \textbf{Right}: {\reftwo Error in the mean value $\big|\int_{\hat{\Gamma}} \hat{c}_s^{(k), \text{fine}}(\hat{x},t) - \hat{c}_s^*(\hat{x},t) \, \text{d}\hat{x}\big|$} over time.
\label{fig:parareal_dist}
}
\end{figure}

\begin{table}[t]
\begin{center} \small
\begin{tabular}{|l|rrrrr|r|}
\hline
k &$P=10$ &$P=20$ &$P=30$ &$P=40$ &$P=50$ &ref. (serial) \\
\hline
1 &$2.80\cdot 10^{-3}$ &$1.42\cdot 10^{-3}$ &$7.49\cdot 10^{-4}$ &$6.53\cdot 10^{-4}$ &$5.04\cdot 10^{-4}$ &${\reftwo 0}$\\ 
2 &$8.73\cdot 10^{-5}$ &$9.46\cdot 10^{-5}$ &$1.27\cdot 10^{-4}$ &$6.22\cdot 10^{-5}$ &$4.78\cdot 10^{-5}$&\\
3 &$3.45\cdot 10^{-5}$ &$2.90\cdot 10^{-5}$ &$2.97\cdot 10^{-5}$ &$1.36\cdot 10^{-5}$ &$9.19\cdot 10^{-6}$ &\\
4 &- &$1.93\cdot 10^{-6}$ &$4.45\cdot 10^{-6}$ &- &- & \\
\hline
Est.\ par. &11\,347\,s &9\,661\,s &8\,692\,s &{\bf 6\,914\,s} &7\,491\,s &26\,840\,s\\
speedup &2.4 &2.8 &3.1 &{\bf 3.9} &3.6 &1.0\\
efficiency &{\bf 24}\,\% &14\,\% &10\,\% &10\,\% &7\,\% &100\,\%\\
\hline
\end{tabular}
\end{center}
\caption{\label{tab:cs_dist} Errors $|\hat{c}_s^{(k),{\sfrei\text{fine}}}(T_{\text{end}})-\hat{c}_s^{*}(T_{\text{end}})|$ at the midpoint $\hat{x}_m$ of $\hat{\Gamma}$ for $P=10$ to $50$ for the parareal algorithm (\cref{alg:parareal}) and a serial reference computation applied to the second numerical example (PDE growth model, {\reftwo $\hat{c}_s^{*}(T_{\text{end}})=0.5186632$}).
Estimated parallel runtimes are shown, as well as speedup and efficiency compared to the reference computation (right column). Details on the estimation of the runtimes are given in~\cref{sec:runtimes}. Best numbers are marked in \textbf{bold face}.}
\end{table}

In the following~\cref{subsec:conv} we investigate {\reftwo numerically} the convergence behavior of the standard parareal algorithm (\cref{alg:parareal}) and the re-usage variant (\cref{alg:para}). Then, we discuss the computational costs of the algorithms applied to the PDE growth model in~\cref{sec:theo_cost}. Finally, we compare the algorithms in~\cref{sec:runtimes} based on estimated runtimes of a parallel implementation.

\subsection{Convergence behavior of the parallel time-stepping algorithms} 
\label{subsec:conv}

We consider again the standard parareal algorithm (\cref{alg:parareal}) and the modification with re-usage of growth values (\cref{alg:para}). As stopping criterion, we now choose
\begin{align*}
| \hat{c}_s^{(k+1), {\sfrei \text{fine}}}(\hat{x}_m,T_{\text{end}}) - \hat{c}_s^{(k), {\sfrei \text{fine}}}(\hat{x}_m,T_{\text{end}})| < \epsilon_{\text{par}} = 10^{-4},
\end{align*}
{\reftwo where $\hat{x}_m=(0,-1)$ is the center of the FSI interface $\hat{\Gamma}$}, 
that is, a slightly more strict tolerance compared to the ODE model.

\paragraph{Standard parareal}
In~\cref{fig:parareal_dist}, we show the {\reftwo evolution of the error in} the first three iterates of the parareal algorithm over time for {\reftwo $P=N_c=10$}. More precisely, we plot the {\reftwo error in} the point functional $\hat{c}_s^{{\sfrei (k), \text{fine}}}(\hat{x}_m,t)$, that is, the foam cell concentration at the center of the FSI interface, and the average $\int_{\hat{\Gamma}} \hat{c}_s^{\sfrei (k), \text{fine}}(\hat{x},t)\, \text{d}x$ of the foam cell concentration over time. 
As in the first numerical example, we observe fast convergence towards the reference values. Again, the foam cell concentrations $\hat{c}_s(\hat{x}_m)$ are significantly overestimated in the initialization step, due to the large time-step $\delta T$ in step (I). These are the starting values for the fine problems in the first parareal iteration. 

In~\cref{tab:cs_dist}, we compare the convergence of the function 
values of $\hat{c}_s^{\sfrei (k), \text{fine}}(\hat{x}_m, T_{\text{end}})$ depending on the number of processes $P$. 
The findings are similar to the first numerical example (\cref{tab:csk}). Again, depending on $P$, three to four
 parareal iterations were sufficient to reach the stopping criterion, with a slightly faster convergence behavior for larger $P$.
 
 {\reftwo
 In~\cref{tab.fixedNc}, we show numerical results for a fixed coarse time step $\delta T=\frac{T_{\text{end}}}{40}$, i.e., $N_c=40$ and varying $P\in\{10,20,40\}$. This means that the cost of the coarse propagator (in terms of the number of micro problems to be solved) is equal in all three cases. While the convergence behavior is slightly faster for smaller $P$, the stopping criterion was satisfied after 3 parareal iterations in all cases. As the fine propagator is cheaper for larger $P$, the fastest computation is $P=N_c=40$, with a speed-up of 4.3 in terms of the number of micro problems to be solved.}

\begin{table}[t]\reftwo
\begin{center} \small
\begin{tabular}{|l|rrr|r|}
\hline
k &$P=10$ &$P=20$ &$P=40$ & ref. (serial) \\
\hline
1 &$8.51\cdot 10^{-4}$ &$7.60\cdot 10^{-4}$ &$6.53\cdot 10^{-4}$
 &{\reftwo 0}\\
2 &$7.64\cdot 10^{-6}$ &$2.25\cdot 10^{-5}$ &$6.22\cdot 10^{-5}$ &-\\
3 &$4.36\cdot 10^{-7}$ &$3.66\cdot 10^{-6}$&$1.36\cdot 10^{-5}$  &-\\
\hline 
\hline 
\hline
\# mp &460 &310 &{\bf 235} &1\,000\\
speedup &2.2 &3.2 &{\bf 4.3}&1.0\\
efficiency &{\bf 22}\,\% &16\,\% &11\,\% &100\,\%\\
\hline
\end{tabular}
\end{center}
\caption{\label{tab.fixedNc} \reftwo Errors $|\hat{c}_s^{(k),{\text{fine}}}(T_{\text{end}})-\hat{c}_s^{*}(T_{\text{end}})|$ at the midpoint $\hat{x}_m$ of $\hat{\Gamma}$ for a fixed coarse time step ($N_c=40$), $P\in\{10, 20, 40\}$ and the parareal algorithm (\cref{alg:parareal}) applied to the second numerical example (PDE growth model).
Estimated parallel runtimes are shown, as well as speedup and efficiency compared to the reference computation (right column). Best numbers are marked in \textbf{bold face}.
}
\end{table}



\paragraph{Parareal with re-usage of growth values}

In~\cref{tab:reusage_pde},
we show the results for the 
parareal
algorithm with re-usage of growth values; cf.~\cref{sec:reusage}. For all tested $P$, we need $k_{\text{par}}=5$ iterations to satisfy the stopping criterion.
Similar to the first numerical example, these are 1-2 additional iterations compared to the standard parareal algorithm.


\subsection{Theoretical discussion of the computational cost} 
\label{sec:theo_cost}

While using the ODE growth model, it was obvious that the solution of the growth model and the communication could be neglected, this requires some discussion for the PDE growth model, as the reaction-diffusion equation~\cref{imex} needs to be solved to advance the foam cell concentration $\hat{c}_s$.

\paragraph{Standard parareal}

 In the standard parareal algorithm (\cref{alg:parareal}), a time step of the PDE growth model always follows the solution of a micro problem, where the growth values $\gamma(\sigmat_f^{WS})$ are computed. Thus, the numbers of time steps of the PDE growth model and micro problems to be solved coincide. Considering that the latter consists of $\geq 100$ time steps, each involving the solution of a nonlinear, coupled FSI problem,
 while the growth model only requires the solution of a single 
scalar PDE of reaction-diffusion type,
it is clear that the computational cost of the growth model is negligible.  

Concerning communication, each process $p=1,...,P$ needs to {\reftwo communicate} its final foam cell concentration ${\cal F}(\overline{c}_s^{(k)}(T_p)) = \hat{c}_s^{\sfrei (k), \text{fine}}(t_{p,n_p})$, which is a scalar-valued finite element function defined in the solid domain $\hat{\SO}$, {\reftwo to the master process if the master-slave approach is used. No communication of ${\cal F}(\overline{c}_s^k(T_p))$ is required in a distributed approach; cf.~\cref{sec:parallelization}.}
Then, after the coarse propagator, 
the variables $\overline{c}_s^{(k+1)}(t_{p\cdot n_p})$ {\reftwo have to be communicated: in the master-slave approach, the data is transferred back to the slaves, and in the distributed case, it is communicated to the next process in line.}
Using the discretization outlined above, this corresponds to 369 degrees of freedom that need to be {\reftwo communicated each time.}

{\reftwo As already discussed in~\cref{sec:parallelization}, in the master-slave case, this communication step is unfavorable because it involves all-to-one and one-to-all communication, whereas the communication pattern for the distributed parallelization only involves neighbor communication, which is more beneficial.}
{\reftwo However, in both cases, the two communication steps are only performed once during each parareal iteration. 
In order to investigate if the communication cost is still negligible for the PDE growth model, parallel numerical experiments based on an actual parallel implementation are needed, especially for realistic three-dimensional problems; see also the discussion in~\cite{gotschel_twelve_2021}. Due to the large computational cost of the micro problems in realistic three-dimensional problems, the communication cost might still be comparably small, which motivates the focus on the number of micro problems to be solved in our discussion; recall that, in our case, the micro problems consist of $\geq 100$ time steps of a fully coupled FSI problem.
Obviously, this argumentation needs to be confirmed based on an actual parallel implementation; we will investigate this further in our future work.}

\paragraph{Parareal with re-usage of growth values} The objective of the re-usage algorithm is to reduce the number of micro problems to be solved in the coarse-scale propagator. As this allows for a smaller time-step size $\delta t$ in the coarse problem, the number of reaction-diffusion equations 
to be solved, might increase significantly. Instead of {\reftwo $N_c$ of} such equations in step (II.b)(i) of~\cref{alg:parareal},~\cref{alg:para} requires $N_f$ reaction-diffusion steps in step (II.b)(ii), where $N_f=1\,000$ in the example used in this section. Thus, in $k_{\text{par}}$ iterations of the parareal algorithm, the number of such equations to be solved is $k_{\text{par}}\cdot N_f$ in step (II.b)(ii) plus {\reftwo $N_c$} equations in step (I) and $k_{\text{par}} \cdot \lceil \frac{N_f}{P}\rceil$ in step (II.a)(ii) of~\cref{alg:para}. In total,~\cref{alg:para} requires the solution of 
\begin{align*}
k_{\text{par}}\cdot \left(N_f + \lceil N_f/P\rceil\right) + {\reftwo N_c}
\end{align*} 
reaction-diffusion equations. For a comparison, we note that the number of micro problems to be solved for $P$ processes was $(k_{\text{par}}+1) {\reftwo N_c} + k_{\text{par}}\cdot \lceil N_f /P \rceil$ (see~\cref{sec:compcost}). This means that the number of reaction-diffusion equations to be solved is by a factor
\begin{align*}
\frac{k_{\text{par}}\cdot \left(N_f + \lceil N_f/P\rceil\right) + {\reftwo N_c}}{(k_{\text{par}}+1) {\reftwo N_c} + k_{\text{par}}\cdot \lceil N_f /P \rceil}
\end{align*}
 larger compared to the number of micro-problems. A simple calculation yields that
\begin{align*}
\frac{k_{\text{par}}\cdot \left(N_f + \lceil N_f/P\rceil\right) + {\reftwo N_c}}{(k_{\text{par}}+1) {\reftwo N_c} + k_{\text{par}}\cdot \lceil N_f /P \rceil} \,\leq\, \frac{k_{\text{par}} \left(N_f + \lceil N_f/P\rceil + {\reftwo N_c}\right)}{k_{\text{par}} \left({\reftwo N_c} + \lceil N_f /P \rceil\right)} \,=\,  \frac{N_f}{{\reftwo N_c}+\lceil N_f /P \rceil} + 1 \,\leq\,\frac{\sqrt{N_f}}{2} +1.
\end{align*} 
{\reftwo In the last inequality, we have used that $N_c \geq P$.}
For $N_f=1\,000$, the bound on the right-hand side is approximately $16.8$. 

Of course, another option would be to use coarser time steps for the coarse propagator. However,
noting again that a micro problem consists of $\geq 100$ FSI steps, the computational cost for $\leq 16.8$ scalar reaction-diffusion equations is still much cheaper. This will be confirmed in the following section, where we show estimated runtimes of a parallel implementation and the respective contributions from the micro problems and the solves of reaction-diffusion equations.

{\reftwo For the re-usage of growth values and the distributed parallelization scheme, the communication cost does not change compared to the standard parareal algorithm. This is because the growth values to be re-used are already available on the process; this follows directly from the discussion in~\cref{sec:reusage}. In the case, of master-slave communication, the communication cost is increased. In particular,}
$N_f$ growth functions $\overline{\gamma}_i, i=1,...,N_f$ need to be transferred from the processes $p=1,...,P$ to the master process. 
Each $\overline{\gamma}_i$ is a spatially discretized function which is defined on the FSI interface $\hat{\Gamma}$. In the example considered here, it is non-zero only in the area $\hat{x}_1\in [-1,1]$ which corresponds to 7 degrees of freedom in our discretization; again, for a realistic three-dimensional problem, the number of degrees of freedom will increase drastically but will still remain small compared to the full problem size. {\reftwo As mentioned before, we assume that the communication cost is negligible compared to the solution of the micro problems in our discussion. Again, this assumption has to be tested in the future based on parallel results.}

\begin{table}[t]
\begin{center} \small
\hspace*{-0.5cm}
\setlength{\tabcolsep}{4pt}\renewcommand{\arraystretch}{1.1}
\begin{tabular}{|l|rrrrrrr|r|}
\hline
k &$P=10$ &$P=20$ &$P=30$ &$P=40$ &$P=50$ &$P=60$ &$P=70$ &ref. \\
\hline
1 &\small $2.80\cdot 10^{-3}$ &\small $1.42\cdot 10^{-3}$ &\small $7.49\cdot 10^{-4}$ &\small $6.53\cdot 10^{-4}$ &\small $5.04\cdot 10^{-4}$ &\small $4.00\cdot 10^{-4}$ &\small $3.32\cdot 10^{-4}$  &\small {\reftwo 0}\\
2 &\small $6.37\cdot 10^{-4}$ &\small $6.96\cdot 10^{-4}$ &\small $5.62\cdot 10^{-4}$ &\small $4.08\cdot 10^{-4}$ &\small $3.45\cdot 10^{-4}$ &\small $3.09\cdot 10^{-4}$ &\small $2.82\cdot 10^{-4}$  & \\ 
3 &\small $1.33\cdot 10^{-4}$ &\small $1.87\cdot 10^{-4}$ &\small $2.06\cdot 10^{-4}$ &\small $1.56\cdot 10^{-4}$ &\small $1.39\cdot 10^{-4}$ &\small $1.30\cdot 10^{-4}$ &\small $1.25\cdot 10^{-4}$ & \\ 
4 &\small $2.87\cdot 10^{-5}$ &\small $5.31\cdot 10^{-5}$ &\small $6.14\cdot 10^{-5}$ &\small $5.08\cdot 10^{-5}$ &\small $4.61\cdot 10^{-5}$ &\small $4.29\cdot 10^{-5}$ &\small $4.02\cdot 10^{-5}$ & \\ 
5 &\small $5.05\cdot 10^{-6}$ &\small $1.47\cdot 10^{-5}$ &\small $1.67\cdot 10^{-5}$ &\small $1.45\cdot 10^{-5}$ &\small $1.20\cdot 10^{-5}$ &\small $1.07\cdot 10^{-5}$ &\small $1.05\cdot 10^{-5}$  &\\
\hline
Est.\ par. &17\,733\,s &9\,685\,s &7\,105\,s &5\,902\,s &5\,277\,s &4\,925\,s &{\bf 4\,804\,s} &26\,840\,s\\
speedup      &1.5     &2.8    &3.8    &4.5 &5.1 &5.4 &{\bf 5.6} &1.0\\
efficiency   &{\bf 15}\,\%  &14\,\%  &13\,\% &11\,\% &10\,\% &9\,\% &8\,\% &100\,\%\\
\hline
\end{tabular}
\caption{\label{tab:reusage_pde} Errors $|c_s^{(k), {\sfrei \text{fine}}}(\hat{x}_m,T_{\text{end}})-\hat{c}_s^{*}(\hat{x}_m,T_{\text{end}})|$ at the midpoint $\hat{x}_m$ of $\hat{\Gamma}$ for $P=10$ to $70$ for~\cref{alg:para} (Re-usage of growth values) applied to the second numerical example (PDE growth model). We show estimated parallel runtimes, as well as speedup and efficiency compared to the reference computation (right column). Details on the estimation of runtimes are given in~\cref{sec:runtimes}. Best numbers are marked in \textbf{bold face}.
}
\end{center}
\end{table}

\subsection{Comparison of runtimes} \label{sec:runtimes}


The computational results given in this paper serve as a proof of concept 
to test the 
presented algorithms. From the discussion in the previous subsection, {\reftwo we assume that} the only relevant computational cost comes from the solution of the micro problems. 
As mentioned before 
in our current implementation, we do not solve the fine-scale problems in parallel on different processes. Instead, they are solved
sequentially one after the other on a single 
process. Since the cost for the micro problems dominates the computation times for any serial or parallel simulation of the plaque growth, we can still discuss the parallelization potential of the methods.
The parallel implementation of the algorithms itself is subject to future work.

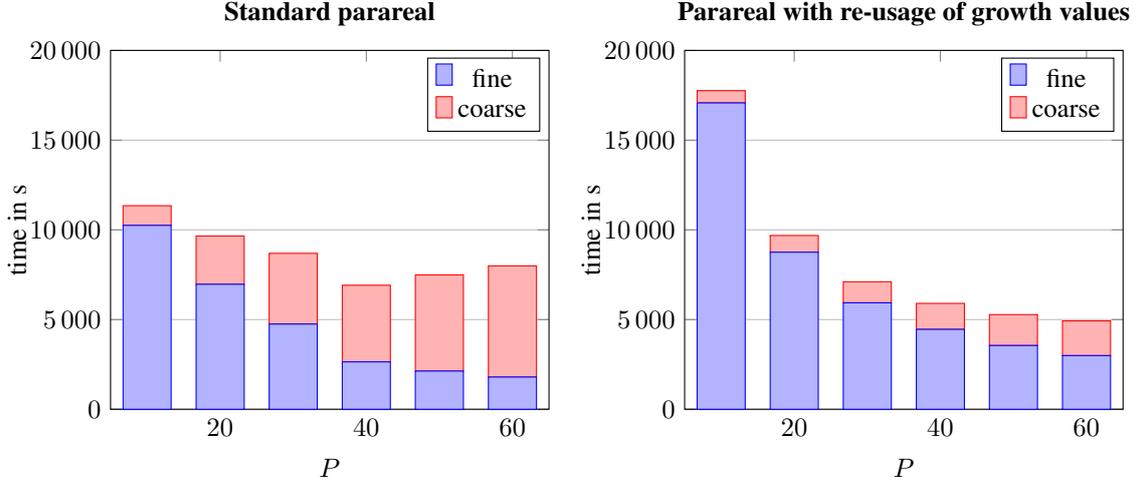
\begin{figure}
\begin{center}
\begin{tikzpicture}
\begin{axis}[
	width=0.49\textwidth,
	title={\textbf{Standard parareal}},
	ybar stacked,
	bar width=18pt,
	xlabel={$P$},
	ymin=0,ymax=20000,
	ylabel={time in s},
	yticklabel style={
        /pgf/number format/fixed,
        /pgf/number format/1000 sep={\,},
        /pgf/number format/precision=5
        },
	scaled y ticks=false,
	ymajorgrids
]
\addplot[color=blue,fill=blue!30] 
	coordinates {(10,10251) (20,6968) (30,4749) (40,2641) (50,2130) (60,1796)};
\addplot[color=red,fill=red!30]
	coordinates {(10,1096) (20,2691) (30,3943) (40,4273) (50,5361) (60,6197)};
\legend{fine, coarse}
\end{axis}
\end{tikzpicture}
\hfill
\begin{tikzpicture}
\begin{axis}[
	width=0.49\textwidth,
	title={\textbf{Parareal with re-usage of growth values}},
	ybar stacked,
	bar width=18pt,
	xlabel={$P$},
	ymin=0,ymax=20000,
	ylabel={time in s},
	yticklabel style={
        /pgf/number format/fixed,
        /pgf/number format/1000 sep={\,},
        /pgf/number format/precision=5
        },
	scaled y ticks=false,
	ymajorgrids
]
\addplot[color=blue,fill=blue!30] 
	coordinates {(10,17075) (20,8755) (30,5934) (40,4454) (50,3555) (60,2992)};
\addplot[color=red,fill=red!30] 
	coordinates {(10,685) (20,930) (30,1171) (40,1448) (50,1722) (60,1933)};
\legend{fine, coarse}
\end{axis}
\end{tikzpicture}
\end{center}
\caption{Illustration of the computational times spent within the coarse- ({\reftwo in serial}) and fine-scale ({\reftwo in parallel}) problems in~\cref{alg:parareal} and~\cref{alg:para} for the second numerical example (PDE growth problem)
\label{fig:runtimes_master_slave} 
}
\end{figure}

\begin{table}[t]
\begin{center} 
\begin{minipage}{0.49\textwidth} \centering \small
\textbf{Standard parareal} 
(\cref{alg:parareal})\\[1mm]
\begin{tabular}{|l|rr@{\hspace{5pt}}r|r|}
\hline
P &\textbf{\reftwo coarse} &\textbf{\reftwo fine}: max. & (aver.) & est.\ par. \\ \hline
10 &1\,096\,s &10\,251\,s & (8\,009\,s) &11\,347\,s \\
20 &2\,691\,s &6\,968\,s & (5\,318\,s) &9\,661\,s\\
30 &3\,943\,s &4\,749\,s & (3\,581\,s) &8\,692\,s\\
40 &4\,273\,s &2\,641\,s & (1\,988\,s) &6\,914\,s\\
50 &5\,361\,s &2\,130\,s & (1\,601\,s) &7\,491\,s\\
60 &6\,197\,s &1\,796\,s & (1\,331\,s) &7\,993\,s\\
\hline
\end{tabular}
\end{minipage}
\begin{minipage}{0.49\textwidth} \centering \small
\textbf{Parareal 
with re-usage of growth values} (\cref{alg:para})\\[1mm]
\begin{tabular}{|l|rr@{\hspace{5pt}}r|r|}
\hline
P& \textbf{\reftwo coarse} &\textbf{\reftwo fine}: max. & (aver.) & est.\ par. \\ \hline
10 &\textbf{658\,s} &17\,075\,s & (13\,342\,s) &17\,733\,s \\ 
20 &930\,s &8\,755\,s & (6\,669\,s) &9\,685\,s \\
30 &1\,171\,s &5\,934\,s & (4\,455\,s) &7\,105\,s \\
40 &1\,448\,s &4\,454\,s & (3\,339\,s) &5\,902\,s\\ 
50 &1\,722\,s &3\,555\,s & (2\,670\,s) &5\,277\,s \\
60 &1\,933\,s &2\,992\,s & (2\,219\,s) &4\,925\,s \\
\hline
\end{tabular}
\end{minipage}
\end{center}

\caption{Estimated parallel runtimes (in s) of the parareal algorithm and the variant "Re-usage of growth values" for $P=10,...,60$ and a PDE growth model. We show the time in seconds spent on the master process and the maximum and average time spent on the slave processes. The estimated parallel runtime is the sum of the time spent in serial (master) and the maximum time needed among the slaves;
best numbers marked in \textbf{bold face}. A visualization of the runtimes is given in~\cref{fig:runtimes_mp_pde}.~\label{tab:runtimes_parareal}
}
\end{table}

In this section, 
we will 
give an estimate of the runtimes that would be required in a parallel implementation. For this purpose, we list the serial part and parallel contributions of the computing times in~\cref{tab:runtimes_parareal}. The serial part corresponds to the coarse propagation, which is performed on a single master process (steps (I) and (II.b) in~\cref{alg:parareal} or~\cref{alg:para}). The parallel part corresponds to the solution of the micro problems (step (II.a)). It can be performed concurrently on all processes $p=1,...,P$. The computing times vary slightly across the processes. It increases for larger processes $p$, as more Newton iterations are required for the FSI problem, when the channel is already significantly narrowed due to the advancing plaque growth. The estimated parallel runtimes given in~\cref{tab:cs_dist,tab:reusage_pde,tab:runtimes_parareal} are the sum of the serial part, which corresponds to the coarse propagator, and the maximum time needed by one of the processes $p=1,...,P$ in the parallel part.





\begin{figure}
\begin{center}
\begin{tikzpicture}
\begin{axis}[
	width=0.49\textwidth,
	title={\textbf{Standard parareal}},
	ybar stacked,
	bar width=18pt,
	xlabel={$P$},
	ymin=0,ymax=20000,
	ylabel={time in s},
	yticklabel style={
        /pgf/number format/fixed,
        /pgf/number format/1000 sep={\,},
        /pgf/number format/precision=5
        },
	scaled y ticks=false,
	ymajorgrids
]
\addplot[color=green!50!black,fill=green!50!black!30]
	coordinates {(10,10227) (20,6952) (30,4739) (40,2635) (50,2126) (60,1792)};
\addplot[color=green!80!black,fill=green!80!black!30]
	coordinates {(10,1093) (20,2683) (30,3931) (40,4260) (50,5345) (60,6178)};
\addplot[color=yellow!30!black,fill=yellow!30!black!30]
	coordinates {(10,22.6) (20,14.9) (30,10.4) (40,5.6) (50,4.5) (60,3.8)};
\addplot[color=yellow!80!black,fill=yellow!80!black!30]
	coordinates {(10,3.0) (20,7.4) (30,11.2) (40,11.9) (50,15.0) (60,17.6)};
\legend{mp (fine), mp (coarse), PDE growth (fine), PDE growth (coarse)}
\end{axis}
\end{tikzpicture}
\hfill
\begin{tikzpicture}
\begin{axis}[
	width=0.49\textwidth,
	title={\textbf{Parareal with re-usage of growth values}},
	ybar stacked,
	bar width=18pt,
	xlabel={$P$},
	ymin=0,ymax=20000,
	ylabel={time in s},
	yticklabel style={
        /pgf/number format/fixed,
        /pgf/number format/1000 sep={\,},
        /pgf/number format/precision=5
        },
	scaled y ticks=false,
	ymajorgrids
]
\addplot[color=green!50!black,fill=green!50!black!30] 
	coordinates {(10,17036) (20,8735) (30,5921) (40,4445) (50,3547) (60,2985)};
\addplot[color=green!80!black,fill=green!80!black!30]
	coordinates {(10,284) (20,550) (30,794) (40,1067) (50,1342) (60,1555)};
\addplot[color=yellow!30!black,fill=yellow!30!black!30]
	coordinates {(10,37.4) (20,18.8) (30,12.7) (40,9.5) (50,7.5) (60,6.4)};
\addplot[color=yellow!80!black,fill=yellow!80!black!30]
	coordinates {(10,368.1) (20,374.2) (30,371.2) (40,375.3) (50,374.0) (60,371.8)};
\legend{mp (fine), mp (coarse), PDE growth (fine), PDE growth (coarse)}
\end{axis}
\end{tikzpicture}
\end{center}
\caption{Illustration of the computational times spent within micro-scale problems and PDE growth problems on the coarse ({\reftwo in serial}) and fine scale ({\reftwo in parallel}) in~\cref{alg:parareal} and~\cref{alg:para} for the second numerical example (PDE growth problem). Note that the time needed for the PDE growth model within the fine-scale problems is so small that it is not visible in both plots. The corresponding times can be found in~\cref{tab:runtimes_parareal}.
\label{fig:runtimes_mp_pde}}
\end{figure}
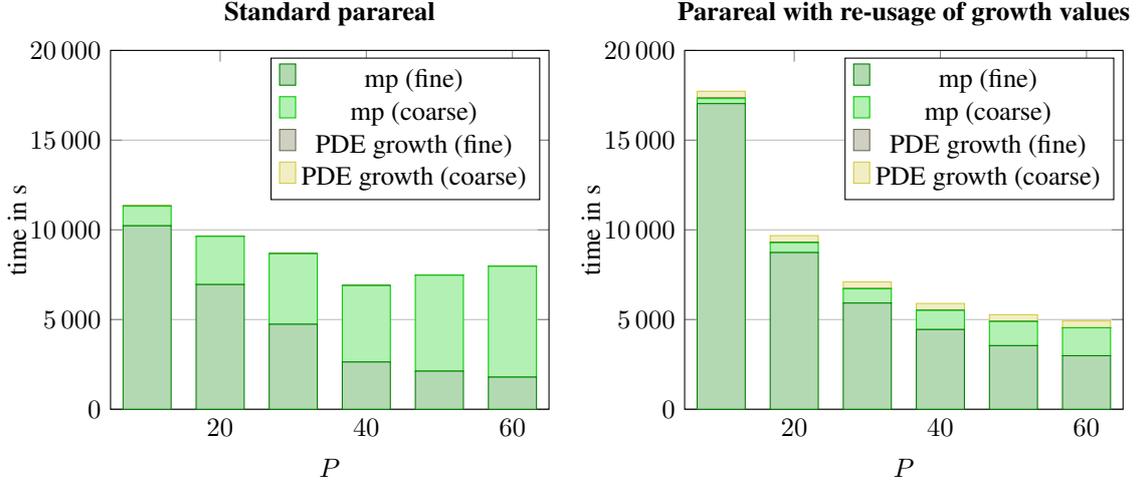

In~\cref{fig:runtimes_master_slave}, we illustrate the computational times needed to solve the fine problems (in parallel) and coarse-scale propagators (in serial) for different processes and with~\cref{alg:parareal} and~\cref{alg:para}. We see that the total computational times decrease for increasing $P$ until $P=40$ for standard parareal and until $P=60$ for the re-usage variant.
In the standard parareal algorithm (\cref{alg:parareal}) the cost of the coarse propagators becomes significant already from $P=20$, while it is much smaller for the re-usage variant. While for $P=10$,~\cref{alg:parareal} is still faster, this changes for $P> 20$. 


The corresponding speed-ups and efficiencies compared to a serial time-stepping are given in the last rows of~\cref{tab:cs_dist,tab:reusage_pde} for standard parareal algorithm~(\cref{alg:parareal}) and the re-usage variant~(\cref{alg:para}), respectively. The lowest computing time for standard parareal is 6\,914\,s for $P=40$, which corresponds to a speed-up of 3.9. For the re-usage variant, we obtain a maximum speed-up of 5.6 for $P=40$, which corresponds to an estimated parallel runtime of 4\,804\,s. The speed-ups are slightly lower than the speed-ups in~\cref{sec:parallel} (cf.~\cref{tab:csk,tab:para2}), that were computed based on the number of micro problems to be solved. This is mostly due to the load imbalancing among the slave processes, as - depending on the state of the plaque growth - some micro problems are more costly than others; for instance, due to a higher number of Newton iterations. The differences can be inferred best by comparing the average and the maximum time spent on the slave processes in~\cref{tab:runtimes_parareal}. 
The efficiencies decrease again monotonically for increasing $P$ in both \cref{tab:cs_dist,tab:reusage_pde}. For the re-usage variant, the efficiencies are again much more stable, due to the reduced cost of the coarse propagator.


Finally, we compare in~\cref{fig:runtimes_mp_pde} the computing times needed for the solution of the FSI micro problems and those needed for the reaction-diffusion equations, for both algorithms and $P=10,...,60$. We observe that the times needed for the latter are negligible in all cases, which confirms the discussion in~\cref{sec:theo_cost}.

\section{Conclusion} \label{sec:conclusion}

We have derived a parareal algorithm for the time parallelization of the macro scale in a two-scale formulation for the simulation of 
atherosclerotic plaque growth. To 
reduce the computational cost of the coarse-scale propagators, we have introduced a variant which re-uses growth values that were computed within the fine problems and avoids 
additional
costly micro-scale computations in serial.

The approaches have been tested on
two different numerical examples of increasing complexity: first, by means of a simple ODE growth model, and secondly, by a PDE model of reaction-diffusion type. In this proof-of-concept, we analyze the approaches based on results and timings of a serial implementation. Since the number of communication steps is low compared to the computational work, we are still able to provide meaningful results. 
In the first case, we obtain estimated speedups up to 6.3 (standard parareal) respectively 7.8 (re-usage variant) in terms of the number of micro problems to be solved. In the PDE model, the maximum estimated speedups, now based on estimated parallel runtimes, were 3.9 respectively 5.6. 
In future work, these results will have to confirmed using a parallel implementation. 
{\reftwo Therefore, we have discussed theoretically two parallelization schemes, master-slave and distributed parallelization. The discussion indicated that the distributed scheme might be beneficial, which also needs to be tested in future work.}

Additional research should also be invested into
further 
improving the efficiency
of the coarse propagator.  
Therefore, the idea
of~\cref{alg:para} could be extended, for example, by storing and interpolating the computed growth values, instead of simply re-using them. 
{\reftwo For the ODE case, such an interpolation approach has already been applied in our proceedings paper~\cite{FreiHeinleinProceedings} with promising results. The extension to the PDE model is, however, not straight-forward, as an operator ${\cal M}: c_s \to \overline{\gamma}(\sigma_f^{WS}(c_s))$ between spatially distributed functions needs to be approximated.}
Further 
future investigations
include the application of the algorithms in complex three-dimensional geometries, more complex plaque growth, and arterial wall models. Finally, the approaches presented here can be combined with a spatial parallelization of the FSI problems or with adaptive time-stepping on micro and macro scale, as presented by Richter \& Lautsch~\cite{RichterLautsch2021}.

\bibliographystyle{plain}


{\reftwo \section*{Appendix}

In this section, we show convergence of the re-usage parareal algorithm, \cref{alg:para} resp.~\cref{eq:parareal_app}, applied to the ODE model in~\cref{eq:EulerODEcont} 
and its approximation by the explicit Euler method in~\cref{eq:EulerODE2}.
First, we obtain from~\cref{eq:EulerODE2,Lipschitz} directly for $c_s^1, c_s^2 \in \mathbb{R}$
\begin{align}
|{\cal C}(I_n,c_s^1(T_{n-1}), \gamma) - {\cal C}(I_n,c_s^2(T_{n-1}), \gamma)| &= |c_s^1(T_{n-1}) - c_s^2(T_{n-1})|, \label{Lip1}\\
|{\cal C}(I_n,c_s, \gamma(c_s^1)) - {\cal C}(I_n,c_s, \gamma(c_s^2))| &= \delta T |\gamma(c_s^1) - \gamma(c_s^2)| \leq L \delta T |c_s^1 - c_s^2| \label{Lip2}. 
\end{align}
By $c_{s,n-1}(T)$ we denote in the following the function that solves the (continuous) ODE~\cref{eq:EulerODEcont} with initial value $c_{s,n-1}(T_{n-1})=c_{s,n-1}$. We have using~\cref{eq:EulerODE2,Lipschitz,exact}
\begin{align}
\begin{split}\label{long1}
&\Big|\left({\cal F}(I_n,c_s(T_{n-1}), \gamma(c_s(T_{n-1}))) - {\cal C}(I_n,c_s(T_{n-1}), \gamma(c_s(T_{n-1})))\right)\\ 
&\qquad-\left( {\cal F}(I_n,c_{s,n-1}, \gamma(c_{s,n-1})) - {\cal C}(I_n,c_{s,n-1}, \gamma(c_{s,n-1}))\right)\Big|\\
&= \Big|\int_{T_{n-1}}^{T_n} \gamma(c_s(T))) - \gamma(c_{s,n-1}(T)) \, \text{d}T - \delta T \left( \gamma(c_s(T_{n-1})) - \gamma(c_{s,n-1})\right)\Big|\\
&\leq c \delta T^2 \left|{\frac{d}{dt}} \gamma(c_s(T_{n-1}) - {\frac{d}{dt}} \gamma(c_{s,n-1})\right| 
\leq \alpha_0 \delta T^2 \left|c_s(T_{n-1}) - c_{s,n-1}\right|
\end{split}
\end{align}
for some constant $\alpha_0=cL>0$; see also~\cite{GanderHairer}.

The following recursion builds the basis for the error estimate:
\begin{lemma}\label{lem.rek}
Let $c_s^{*}\in C^1(0,T_\text{end})$ be the exact solution of~\cref{eq:EulerODEcont}. $\{\overline{c}_s^{(k)}(T_n)\}_{n=1}^{P}$ be the $k$-th iterate of the re-usage algorithm in~\cref{eq:parareal_app} using the forward Euler method in~\cref{eq:EulerODE2} and let $e_n^{(k)} = |c_s^{*}(T_n) - \overline{c}_{s,n}^{(k)}|$ be the error in the $k$-the iteration of the re-usage algorithm. Under the assumptions made in~\cref{subsec.conv}, it holds for $k=1$ that
\begin{align}\label{keq1}
e_n^{(1)} &\leq \alpha_0 \delta T^2 e_{n-1}^{(0)} +\alpha_1\delta T e_{n-2}^{(0)} + e_{n-1}^{(1)},
\end{align}
and for $k\geq 2$,
\begin{align}\label{eq:kgeq2}
e_n^{(k)} &\leq \alpha_1 \delta T \left(e_{n-1}^{(k-1)} + e_{n-2}^{(k-1)} + e_{n-2}^{(k-2)}\right) + e_{n-1}^{(k)},
\end{align}
with a constant $\alpha_1 \geq \max\{L + \alpha_0 \delta T, L (1+L\delta T)\}$.
\end{lemma}
\begin{proof}
Consider first the case $k\geq 2$. Using~\cref{eq:parareal_app,exact}, we have 
\begin{align*}
e_n^{(k)} &= \big|{\cal F}(I_n,c_s^*(T_{n-1}), \gamma(c_s^*(T_{n-1})) - {\cal C} (I_n,\overline{c}_{s,n-1}^{(k)}, \gamma(c_{s,n-1}^{(k), \text{fine}}))\\
&\qquad-{\cal F}(I_n, \overline{c}_{s,n-1}^{(k-1)}, \gamma(\overline{c}_{s,n-1}^{(k-1)})) + {\cal C}(I_n,\overline{c}_{s,n-1}^{(k-1)}, \gamma(c_{s,n-1}^{(k-1),\text{fine}})) \big|\\
&\leq \big|{\cal F}(I_n,c_s^*(T_{n-1}), \gamma(c_s^*(T_{n-1}))) -
{\cal C} (I_n,c_s^*(T_{n-1}), \gamma(c_s^*(T_{n-1})))\\
&\qquad- {\cal F}(I_n, \overline{c}_{s,n-1}^{(k-1)}, \gamma(\overline{c}_{s,n-1}^{(k-1)})) 
+{\cal C}(I_n, \overline{c}_{s,n-1}^{(k-1)}, \gamma(\overline{c}_{s,n-1}^{(k-1)})) \big|
\\
&\quad+ \big| {\cal C}(I_n, \overline{c}_{s,n-1}^{(k-1)}, \gamma(\overline{c}_{s,n-1}^{(k-1)}))-{\cal C}(I_n,\overline{c}_{s,n-1}^{(k-1)}, \gamma(c_{s,n-1}^{(k-1),\text{fine}})) \big| \\
&\quad+\big| {\cal C} (I_n,\overline{c}_{s,n-1}^{(k)}, \gamma(c_{s,n-1}^{(k), \text{fine}}))
- {\cal C} (I_n,\overline{c}_{s,n-1}^{(k)}, \gamma(c_s^*(T_{n-1}))) \big|\\
&\quad+ \big| {\cal C} (I_n,\overline{c}_{s,n-1}^{(k)}, \gamma(c_s^*(T_{n-1}))) -
{\cal C} (I_n,c_s^*(T_{n-1}), \gamma(c_s^*(T_{n-1})))
\big|
\end{align*}
From~\cref{Lip1,Lip2,long1}, we have
\begin{align}\label{eq:almostfinal}
\begin{split}
e_n^{(k)} &\leq \alpha_0 \delta T^2 \left|c_s^*(T_{n-1}) - \overline{c}_{s,n-1}^{(k-1)}\right|
+L \delta T \big|\overline{c}_{s,n-1}^{(k-1)} - c_{s,n-1}^{(k-1),\text{fine}}\big| \\
&\quad + L\delta T\big|c_{s,n-1}^{(k), \text{fine}} -    c_s^*(T_{n-1}) \big| +\big| \overline{c}_{s,n-1}^{(k)} - c_s^*(T_{n-1})\big| \\
&= \alpha_0 \delta T^2 e_{n-1}^{(k-1)} + L \delta T ( e_{n-1}^{(k-1),\text{fine}} + e_{n-1}^{(k-1)})+ L\delta T e_{n-1}^{(k),\text{fine}} + e_{n-1}^{(k)}, 
\end{split}
\end{align} 
where $e_{n-1}^{(k),\text{fine}} := | c_{s,n-1}^{(k), \text{fine}} - c_s^*(T_n)|$. 
This error contribution is estimated further by using~\cref{Defcsfine,exact,Lip1}
\begin{align}
\begin{split}\label{eq:finearg}
e_{n-1}^{(k),\text{fine}} := | c_{s,n-1}^{(k), \text{fine}} - c_s^*(T_n)| &=| {\cal F}(I_{n-1}, c_{s,n-2}^{(k-1)}, \gamma(c_{s,n-2}^{(k-1)})) - {\cal F}(I_{n-1}, c_s^*(T_{n-2}), \gamma(c_s^*(T_{n-2}))|\\
&=\big| c_{s,n-2}^{(k-1)} - c_s^*(T_{n-2}) + \int_{T_{n-2}}^{T_{n-1}} \gamma(c_{s,n-2}^{(k-1)}(t))) - \gamma(c_s^*(t))\,\text{d}t\big|\\
&\leq (1+L\delta T) e_{n-2}^{(k-1)}.
\end{split}
\end{align} 
Inserting this into~\cref{eq:almostfinal} yields~\cref{eq:kgeq2}. For $k=1$, we have $c_{s,n-1}^{(0),\text{fine}} = c_{s,n-1}^{(0)}$. An analogous argumentation yields
\begin{align*}
e_n^{(1)} &= \big|{\cal F}(I_n,c_s^*(T_{n-1}), \gamma(c_s^*(T_{n-1}))) - {\cal C} (I_n,\overline{c}_{s,n-1}^{(1)}, \gamma(c_{s,n-1}^{(1),\text{fine}}))\\
&\qquad-{\cal F}(I_n, \overline{c}_{s,n-1}^{(0)}, \gamma(\overline{c}_{s,n-1}^{(0)})) + {\cal C}(I_n,\overline{c}_{s,n-1}^{(0)}, \gamma(c_{s,n-1}^{(0)})) \big|\\
&\leq \big|{\cal F}(I_n,c_s^*(T_{n-1}), \gamma(c_s^*(T_{n-1}))) -
{\cal C} (I_n,c_s^*(T_{n-1}), \gamma(c_s^*(T_{n-1})))\\
&\qquad- {\cal F}(I_n, \overline{c}_{s,n-1}^{(0)}, \gamma(\overline{c}_{s,n-1}^{(0)})) 
+{\cal C}(I_n, \overline{c}_{s,n-1}^{(0)}, \gamma(\overline{c}_{s,n-1}^{(0)})) \big|
\\
&\quad+\big| {\cal C} (I_n,\overline{c}_{s,n-1}^{(1)}, \gamma(c_{s,n-1}^{(1), \text{fine}}))
- {\cal C} (I_n,\overline{c}_{s,n-1}^{(1)}, \gamma(c_s^*(T_{n-1}))) \big|\\
&\quad+ \big| {\cal C} (I_n,\overline{c}_{s,n-1}^{(1)}, \gamma(c_s^*(T_{n-1}))) -
{\cal C} (I_n,c_s^*(T_{n-1}), \gamma(c_s^*(T_{n-1})))
\big|\\
&\leq  \alpha_0 \delta T^2 e_{n-1}^{(0)} + L\delta T e_{n-1}^{(1),\text{fine}} + e_{n-1}^{(1)}.
\end{align*}
Then,~\cref{keq1} follows by using~\cref{eq:finearg}.
\end{proof}

\begin{lemma}\label{lem.conv_reusage}

Let $e_n^{(k)} = |\overline{c}_s^{(k)}(T_n)-c_s^{*}(T_n)|$ be the error in the $k$-the iteration of the re-usage algorithm. Under the assumptions made above, it holds for $k\in\mathbb{N}_0$ and $n\in\mathbb{N}$ that
\begin{align}
e_n^{(k)} &\leq L \tilde{\alpha}^{k} \beta^{n-k} 
\cdot \left(\sum_{l=\lceil \frac{k}{2}+1\rceil}^{k+1} \frac{3^{l-1}}{l!} \delta T^{l+1} \left( n-\frac{k}{2}\right)^{l}\right) \max_{t\in [0,T_\text{end}]} |\partial_t c_s^{*}(t)|\label{eq:est1}\\
&\leq L \tilde{\alpha}^{k} \beta^{n-k} \delta T  \max\left\{1,T_{n-\frac{k}{2}}\right\}^k \frac{3^{\lceil \frac{k}{2}\rceil}}{\lceil \nicefrac{k}{2} \rceil !} \max_{t\in [0,T_\text{end}]} |\partial_t c_s^{*}(t)|\label{eq:est2},
\end{align}
where $\tilde{\alpha} = \max\{\alpha_0 \delta T, \alpha_1, 1\}$ and $\beta=1+L\delta T$. 
\end{lemma}

\begin{proof}
We prove the lemma by induction over $k\in\mathbb{N}$.
For $k=0$, a standard estimate of the forward Euler methods gives, using~\cref{Lipschitz},
\begin{align}
\begin{split}\label{eq:en0}
e_n^{(0)} := | c_{s,n}^0 - c_s^{*}(T_n)|
&\leq | c_{s,n-1}^0 - c_s^{*}(T_{n-1})| + \int_{T_{n-1}}^{T_n} |\gamma(c_{s,n-1}^0) - \gamma(c_s^{*}(t))|\, \text{d}t\\
&\leq e_{n-1}^{(0)} + L  \int_{T_{n-1}}^{T_n} |c_{s,n-1}^0 - \underbrace{c_s^{*}(t)}_{=c_s^*(T_{n-1}) + \delta T \partial_t c_s^{*}(\xi)|}| \, \text{d}t\\
&\leq (\underbrace{1+L\delta T}_{=\beta}) e_{n-1}^{(0)} + L  \delta T^2\max_{t\in [0,T_{\text{end}}]} |\partial_t c_s^*(t)|.
\end{split}
\end{align}
To abbreviate the notation, we set $c_0 := \max_{t\in [0,T_\text{end}]} |\partial_t c_s^*(t)|$. We apply~\cref{eq:en0} recursively
\begin{align}\label{eq:keq0}
e_n^{(0)} \leq \beta e_{n-1}^{(0)} + c_0 L \delta T^2
\leq c_0 L \delta T^2 \sum_{l=0}^{n-1} \beta^{l}
=  c_0 L \delta T^2 \frac{\beta^{n}-1}{\underbrace{\beta-1}_{=L\delta T}} &=c_0 \delta T (\beta^n-1)
\leq c_0 L n \delta T^2 \beta^n.
\end{align}
In the last inequality, we have used that $\beta^n-1\leq n L\delta T \beta^n$, which can be shown by induction over $n$. The estimate in~\cref{eq:keq0} is exactly~\cref{eq:est1} for $k=0$.

 For $k=1$, we have using~\cref{keq1,eq:keq0}, the fact that $e_k^{(k)}=0$, and by the definition of the parareal algorithm
 \begin{align*}
e_n^{(1)} &\leq \alpha_0 \delta T^2 e_{n-1}^{(0)} +\alpha_1\delta T e_{n-2}^{(0)} + e_{n-1}^{(1)}\\
&\leq \alpha_0 c_0 L \delta T^4 (n-1) \beta^{n-1}
+\alpha_1 c_0 L \delta T^3 (n-2) \beta^{n-2} + e_{n-1}^{(1)}\\
&\leq 2\tilde{\alpha} c_0 L \delta T^3 (n-1) \beta^{n-1} + e_{n-1}^{(1)}\\
&\leq 2\tilde{\alpha} c_0 L\beta^{n-1} \delta T^3 \sum_{l=1}^{n-1} l + \underbrace{e_{1}^{(1)}}_{=0}\\
&=  2\tilde{\alpha} c_0 L\beta^{n-1} \delta T^3 \frac{n(n-1)}{2}
\leq \tilde{\alpha} c_0 L\beta^{n-1} \delta T^3 \left(n-\nicefrac{1}{2}\right)^2.
 \end{align*}
 This is by a factor $\frac{3}{2}$ smaller compared to~\cref{eq:est1}.


Now, let $k\geq 2$. We assume that the estimate in~\cref{eq:est1} is true for $k-1$ and $k-2$. By~\cref{eq:kgeq2} and the assumption of the induction, we have 
\begin{align}
\begin{split}\label{eq:bigsplit}
e_n^{(k)} &\leq \alpha_1 \delta T \left(e_{n-1}^{(k-1)} + e_{n-2}^{(k-1)} + e_{n-2}^{(k-2)}\right) + e_{n-1}^{(k)}\\
&\leq 
\alpha_1 \delta T \Bigg(c_0 L \tilde{\alpha}^{k-1} \beta^{n-k} 
\cdot \left(\sum_{l=\lceil \frac{k+1}{2}\rceil}^{k} \frac{3^{l-1}}{l!} \delta T^{l+1} \left( n-\frac{k+1}{2}\right)^{l}\right) \\
&\qquad\qquad+ c_0 L \tilde{\alpha}^{k-1} \beta^{n-k-1} 
\cdot \left(\sum_{l=\lceil \frac{k+1}{2}\rceil}^{k} \frac{3^{l-1}}{l!} \delta T^{l+1} \left( n-1-\frac{k+1}{2}\right)^{l}\right)\\
&\qquad\qquad+ c_0 L \tilde{\alpha}^{k-2} \beta^{n-k} 
\cdot \left(\sum_{l=\lceil \frac{k}{2}\rceil}^{k-1} \frac{3^{l-1}}{l!} \delta T^{l+1} \left( n-1-\frac{k}{2}\right)^{l}\right)
\Bigg)\,\,+e_{n-1}^{(k)}\\
&\leq 3 c_0 L \tilde{\alpha}^{k} \beta^{n-k} 
\cdot \left(\sum_{l=\lceil \frac{k}{2}\rceil}^{k} \frac{3^{l-1}}{l!} \delta T^{l+2} \left( n-\frac{k+1}{2}\right)^{l}\right) +e_{n-1}^{(k)}.  
\end{split}
\end{align}
We apply this estimate recursively for the last term $e_{n-1}^{(k)}$ in~\cref{eq:bigsplit} to get
\begin{align*}
e_n^{(k)} &\leq c_0 L \tilde{\alpha}^{k} \beta^{n-k} 
\cdot \sum_{m=k+1}^{n} \left(\sum_{l=\lceil \frac{k}{2}\rceil}^{k} \frac{3^{l}}{l!} \delta T^{l+2} \left( m-\frac{k+1}{2}\right)^{l}\right) +\underbrace{e_{k}^{(k)}}_{=0}\\
&= c_0 L \tilde{\alpha}^{k} \beta^{n-k} 
\cdot \left(\sum_{l=\lceil \frac{k}{2}\rceil}^{k} \frac{3^{l}}{l!} \delta T^{l+2} \cdot \underbrace{\left(\sum_{m=\frac{k+1}{2}}^{n-\frac{k+1}{2}}m^{l}\right)}_{=:s}\right).
\end{align*}
In the case that $\frac{k+1}{2}$ and $n-\frac{k+1}{2}$ are no natural number, the sum $s$ is to be understood in such a way that the index $m$ advances iteratively by 1 until reaching the upper limit.  

The sum $s$ is an approximation of the integral 
\begin{align*}
\int_{\nicefrac{k}{2}}^{n-\nicefrac{k}{2}} x^l\,\text{d}x \leq \frac{1}{l +1} \left(n-\frac{k}{2}\right)^{l+1}
\end{align*}
with the midpoint rule. As $f(x) =  x^l$ is a convex function ($f''(x)>0$) for a positive $x$, the integral is an upper bound for $s$; see the error representation of the midpoint rule, e.g., in Theorem 8.41 in~\cite{RichterWickBook}.

We have thus shown that
\begin{align*}
e_n^{(k)} \leq c_0 L \tilde{\alpha}^{k} \beta^{n-k} 
\cdot \left(\sum_{l=\lceil \frac{k}{2}\rceil}^{k} \frac{3^{l}}{(l+1)!} \delta T^{l+2} \left(n-\frac{k}{2}\right)^{l+1}\right).
\end{align*}
Shifting the index $l$ by 1 gives~\cref{eq:est1}. The estimate in~\cref{eq:est2} follows from~\cref{eq:est1} by noting that 
\begin{align*}
(n-\frac{k}{2}) \delta T  = T_{n-\frac{k}{2}} \leq \max\left\{1,T_{n-\frac{k}{2}}\right\}
\end{align*} and the fact that the term $\frac{3^{l-1}}{l!}$ is decreasing for $l\geq 2$:
\begin{align*}
e_n^{(k)} &\leq c_0 L \tilde{\alpha}^{k} \beta^{n-k} 
\cdot \left(\sum_{l=\lceil \frac{k}{2}+1\rceil}^{k+1} \frac{3^{l-1}}{l!} \delta T^{l+1} \left(n-\frac{k}{2}\right)^{l}\right) \\
&\leq c_0 L \tilde{\alpha}^{k} \beta^{n-k} \delta T \left(\max\{1,T_{n-\frac{k}{2}}\}\right)^{k+1}
\cdot \left(\sum_{l=\lceil \frac{k}{2}+1\rceil}^{k+1} \frac{3^{l-1}}{l!}\right)\\
&\leq c_0 L \tilde{\alpha}^{k} \beta^{n-k} \delta T \left(\max\left\{1,T_{n-\frac{k}{2}}\right\}\right)^{k+1}
\lfloor \frac{k}{2}+1\rfloor \frac{3^{\lceil\frac{k}{2}\rceil}}{\lceil\frac{k}{2}+1\rceil !}\\
&\leq c_0 L \tilde{\alpha}^{k} \beta^{n-k} \delta T \left(\max\left\{1,T_{n-\frac{k}{2}}\right\}\right)^{k+1} \frac{3^{\lceil\frac{k}{2}\rceil}}{\lceil\frac{k}{2}\rceil !}.
\end{align*}

\end{proof}}

\end{document}